\documentclass[10pt,a4paper]{article}

\usepackage[utf8]{inputenc}
\usepackage{amsmath}
\usepackage{amsfonts}
\usepackage{amssymb}
\usepackage{amsthm}
\usepackage[headsep=10mm]{geometry}
\usepackage{caption}
\usepackage{subcaption}
\usepackage[pdftex]{graphicx}
\usepackage{bm}  
\usepackage[hidelinks]{hyperref}
\usepackage{fancyhdr}
\usepackage{stmaryrd}  
\usepackage{enumitem}  
\usepackage{mathtools}
\usepackage{array}  

\usepackage{tikz}
\usetikzlibrary{matrix}

\numberwithin{equation}{section}

\newcommand{\grad}[1]{{\bm{\nabla}{#1}}}
\newcommand{\dv}[1]{\bm{\nabla\cdot{#1}}}
\newcommand{\curl}[1]{\bm{\nabla}\times\bm{#1}}
\newcommand{\pfrac}[2]{\frac{\partial{#1}}{\partial{#2}}}
\newcommand{\dx}[1]{\hspace{0.75mm}\mathrm{d}{#1}}

\newcommand{\cder}[2]{\frac{\mathrm{D}#1}{\mathrm{D}#2}}

\newcommand{\footremember}[2]{%
    \footnote{#2}
    \newcounter{#1}
    \setcounter{#1}{\value{footnote}}%
}
\newcommand{\footrecall}[1]{%
    \footnotemark[\value{#1}]%
}

\begin{document}
\vspace{-1cm}
\title{A Compatible Finite Element Discretisation for the
Moist Compressible Euler Equations}
\author{Thomas M. Bendall\footremember{Imperial}{Mathematics Department, Imperial College London, UK}, Thomas H. Gibson\footrecall{Imperial}, Jemma Shipton\footremember{Exeter}{Mathematics Department, University of Exeter, UK},\\ Colin J. Cotter\footrecall{Imperial}, Ben Shipway\footremember{Met Office}{Dynamics Research, Met Office, Exeter, UK}}
\maketitle
\vspace{-1cm}
\begin{abstract}
\noindent
We present new discretisation of the moist compressible Euler equations, using the compatible finite element framework identified in Cotter and Shipton (2012).
The discretisation strategy is described and details of the parametrisations of moist processes are presented.
A procedure for establishing hydrostatic balance for moist atmospheres is introduced, and the model's performance is demonstrated through several test cases, two of which are new.
\end{abstract}

\section{Introduction and Motivation} \label{sec:intro}
Latitude-longitude grids over the sphere have been popular for use in the dynamical cores of numerical weather prediction models due to the orthogonality of the meridians and circles of latitude.
This orthogonality can be exploited to gain a number of desirable numerical properties \cite{staniforth2012horizontal}.
However, as these grids are refined, the grid points around the poles converge.
On massively-parallel computers, there is a bottleneck in the data communication around the poles, which compromises the scalability of the model.
There has therefore been a search for alternative grids that still provide the numerical properties described in \cite{staniforth2012horizontal}. \\
\\
Compatible finite element methods offer a promising solution to this problem.
These are finite element methods in which the variables lie in different function spaces, such that the discrete equations replicate the vector calculus identities of the continuous equations, such as $\bm{\nabla}\times\grad{f}=\bm{0}$ and $\bm{\nabla\cdot\nabla}\times\bm{f}=0$.
\cite{cotter2012mixed} showed for the linear rotating shallow-water equations that a compatible finite element discretisation, which facilitates the use of non-orthogonal meshes, still maintains many of the properties of \cite{staniforth2012horizontal}.
More recent work using compatible finite element methods for geophysical fluids includes \cite{cotter2014finite}, \cite{natale2016compatible} and \cite{yamazaki2017vertical}.
The UK Met Office is therefore developing a dynamical core using such a discretisation \cite{melvin2019mixed}. \\
\\
The main result of this paper is the presentation of a model solving the three dimensional moist compressible Euler equations within this compatible finite element framework.
The inclusion of moisture introduces to the model a basic form of `physics' (a description of unresolved and diabatic physical processes), which gives a robust test of the compressible finite element discretisation.
We then present the use of this model in a series of test cases featuring moisture.\\
\\
The structure of this paper is as follows.
First we present the form of the continuous equations that we aim to solve in Section \ref{sec:equations}.
Then Section \ref{sec:dynamics} describes the spatial and temporal discretisation of the dynamical part of the model, presenting the function spaces for the variables and the strategy for solving the equations of motion.
Section \ref{sec:parametrisations} then details the specifics of the moisture parametrisations, including a discussion of how to combine fields from different function spaces within a parametrisation.
After presenting our approach to initialising states with hydrostatic balance in Section \ref{sec:hydrostatic balance}, we then demonstrate its use in a number of test cases in Section \ref{sec:test cases}. \\
\\
Finally, the model presented in this paper is written as part of Gusto, a library for dynamical cores using compatible finite element discretisations.
Gusto itself is based on the Firedrake software of \cite{rathgeber2017firedrake}, the development of which is based at Imperial College London.
This software provides automated code generation for finite element methods.
Firedrake also has the functionality for tensor product element and extruded mesh functionality, as described in
\cite{mcrae2016automated}, \cite{bercea2016structure} and \cite{homolya2016parallel}, which we used substantially.
\section{Governing Equations}\label{sec:equations}
We solve the compressible Euler equations featuring three species of moisture: water vapour, cloud water and rain.
The ice phase is neglected.
Motivated by the UK Met Office's most recent dynamical core, ENDGame, \cite{endgame2014documentation}, our prognostic variables are the density of dry air $\rho_d$ and the dry virtual potential temperature $\theta_{vd}$ as prognostic variables, alongside the wind velocity $\bm{v}$ and the mixing ratios $r_v$, $r_c$ and $r_r$.
Here the subscripts respectively denote water vapour, cloud water and rain, whilst the mixing ratio $r_i$ is the ratio of the density by volume of the $i$-th substance to that of dry air, i.e. $r_i:=\rho_i/\rho_d$.
The total mixing ratio of water species is $r_t:=r_v+r_c+r_r$.
The dry potential temperature $\theta_{vd}$ is defined for $r_v$, temperature $T$ and air pressure $p$ by
\begin{equation}
\theta_{vd}:=T\left(\frac{p_R}{p}\right)^\frac{R_d}{c_{pd}}
(1 + r_v/\epsilon)=\theta_d(1+r_v/\epsilon),
\end{equation}
where $c_{pd}$ is the specific heat capacity of dry air, $p_R$ is a reference pressure, $\theta_d$ is the dry potential temperature and $\epsilon:=R_d/R_v$ is the ratio of specific gas constants of dry air and water vapour.
The choice of $\theta_{vd}$ is motivated in
\cite{endgame2014documentation}, which notes that it is the more
natural choice of variable to complement $\rho_d$, and claims
that it may be smoother than the dry potential temperature
$\theta_d$. \\
\\
The full equation set that we use is
\begin{subequations}\label{eqn:wet atmos full}
\begin{align}
 \cder{\bm{v}}{t}+\bm{f}\times\bm{v}+\frac{c_{pd}\theta_{vd}}{1+r_t}\grad{\Pi}+\grad{\Phi}  =& \  \bm{0}, \label{eqn:wet atmos momentum} \\
 \begin{split}
\cder{\theta_{vd}}{t} +\theta_{vd}\left(\frac{R_m}{c_{vml}}-\frac{R_dc_{pml}}{c_{pd}c_{vml}}\right)\dv{v}  =&  \\
-\theta_{vd}\left[\frac{c_{vd}L_v(T)}{c_{vml}c_{pd}T} \right. &
\left. -\frac{R_v}{c_{vml}}\left(1-\frac{R_dc_{vml}}{R_mc_{pd}}\right)-\frac{R_v}{R_m}\right]\cder{r_v}{t}, \label{eqn:wet atmos theta_v eqn of motion} 
\end{split} \\
\cder{\rho_d}{t}+\rho_d\dv{v}  = & \ 0, \\
\cder{r_v}{t}  = & -\dot{r}^c_\mathrm{cond} + 
\dot{r}^r_\mathrm{evap}, \\
\cder{r_c}{t}  = & \ \dot{r}^c_\mathrm{cond} - 
\dot{r}_\mathrm{accr}-\dot{r}_\mathrm{accu}, \\
\cder{r_r}{t} = & \ \dot{r}_\mathrm{accr}-\dot{r}_\mathrm{accu} 
-\dot{r}^r_\mathrm{evap} -S.
\end{align}
\end{subequations}
Here $\Pi$ is the Exner pressure function and $\Phi$ represents the geopotential, whilst $\bm{f}=f\widehat{\bm{k}}$ represents the Coriolis parameter multiplied by the vertical unit vector.
For the specific heat capacities $c_{pd}$, $c_{vd}$, $c_{vml}$ and $c_{pml}$, the specific gas constant $R_m$ and also the latent heat of vaporization $L_v(T)$, we follow closely the values used in \cite{bryan2002benchmark}, which are also listed in the appendix.
The advective derivative is given by
\begin{equation}
\cder{}{t} = \pfrac{}{t}+\bm{v}\bm{\cdot}\grad{}.
\end{equation}
The equation of state can be written as
\begin{equation}
\Pi = \left(\frac{p}{p_R}\right)^\kappa\equiv
\left(\frac{\rho_dR_d\theta_{vd}}{p_R}\right)
^\frac{\kappa}{1-\kappa}, \label{eqn:pi def}
\end{equation}
with $\kappa:=R_d/c_{pd}$. \\
\\
The terms on the left hand sides of (\ref{eqn:wet atmos full})
represent the dynamics, whilst the right hand sides are
considered to be the physics.
The processes $\dot{r}^c_\mathrm{cond}$, $\dot{r}_\mathrm{accr}$,
$\dot{r}^r_\mathrm{evap}$, $\dot{r}_\mathrm{accu}$ and $S$ are
the microphysics parametrisations and are described in
Section \ref{sec:parametrisations}. \\
\\
The momentum equation can also be recast in vector invariant form:
\begin{equation}
\pfrac{\bm{v}}{t}+\left(\curl{v}\right)\times \bm{v} + \bm{f}\times\bm{v} + \frac{1}{2}\grad{v^2}+\frac{c_{pd}\theta_{vd}}{1+r_t}\grad{\Pi}+\grad{\Phi} = \bm{0}.
\end{equation}
The final thing to note is the extra term proportional to
$\dv{v}$ appearing on the left hand side of 
(\ref{eqn:wet atmos theta_v eqn of motion}).
This term is neglected in many models but mentioned in
\cite{thuburn2017use} and \cite{bryan2002benchmark} to be
important in fully capturing convection, particularly in a saturated atmosphere.
In our model it appears in the forcing step of the dynamical 
core.

\section{Dynamics Discretisation}\label{sec:dynamics}
\subsection{Function Spaces}
One of the main results of \cite{cotter2012mixed} was the suggestion of combinations of function spaces for the velocity $\bm{v}$ and the (dry) air density $\rho_d$ within the compatible finite element framework that will satisfy the properties described in \cite{staniforth2012horizontal}.
These spaces mimic the Arakawa C-grid staggering that allowed for steady geostrophic modes, avoidance of spurious pressure modes and more accurate representation of the dispersion relation for Rossby and inertia-gravity waves.
Compatible sets of finite element spaces for $\bm{v}$ and $\rho_d$ will form part of a discrete de Rham complex,  so that $\dv{v}$ is in the same space as $\rho_d$.
Some such families of spaces can be found in the periodic table of finite elements \cite{arnold2014periodic}.\\
\\
In this section we will expand on this to list the function spaces $(V_\rho, V_{\bm{v}}, V_\theta)$ that we use in our model for velocity, density and potential temperature respectively.
We shall consider both two-dimensional vertical slices with quadrilateral cells and also three-dimensional domains using hexahedral elements with quadrilateral faces.
The stratification of geophysical fluids along with their high aspect ratio motivates using an extruded mesh in the vertical, meaning that the mesh has regular layers.
The finite elements on such a mesh can be constructed as tensor product elements, as described in \cite{mcrae2016automated}.
We will also present two configurations of the model, using two different sets of spaces from the same family of finite elements spaces.
We label these two sets as the \textit{lowest-order} $k=0$ spaces and the \textit{next-to-lowest-order} $k=1$ spaces.
More explanation on these can be found in \cite{bendall2019recovered}.
All of these spaces are illustrated in Table \ref{tab:spaces}. \\
\\
The density $\rho_d$ lies in a space that is discontinuous between elements.
For the $k=0$ spaces, density is constant within each cell, while it is linear in each cell for the $k=1$ spaces.
These are the discontinuous Galerkin spaces, denoted by $\mathrm{dQ}_k$.\\
\\
The velocity $\bm{v}$ has continuous normal components between cells but discontinuous tangential components.
The spaces that we use are the Raviart-Thomas spaces, $\mathrm{RT}_k$, \cite{raviart1977mixed}.
For the $k=0$ configuration this means that for each component, the field is linear in the direction of that component and continuous between cells in that direction, but constant within the cell in other directions.
This becomes continuous and quadratic in the direction of the component for the $k=1$ spaces, but discontinuous and linear in other directions.\\
\\
To replicate the Charney-Philips grid, the potential temperature $\theta_{vd}$ is co-located with the vertical component of velocity.
This is the tensor product element of the discontinuous Galerkin element $\mathrm{dQ}_k$ over the horizontal part of the domain with the continuous Galerkin element $\mathrm{Q}_{k+1}$ on the vertical part of the domain.
We therefore denote this space by $\mathrm{dQ}_k\otimes\mathrm{Q}_{k+1}$.
Moisture variables also lie in this space, as they can be considered as a dynamical adjustment to $\theta_{vd}$ and so as to easily facilitate latent heat transfer as the water changes phase.
\\
\begin{table}[h!]
\centering
\begin{tabular}{>{\centering\arraybackslash}m{2cm}|
>{\centering\arraybackslash}m{2.5cm}|
>{\centering\arraybackslash}m{2.5cm}|
>{\centering\arraybackslash}m{2.5cm}|
>{\centering\arraybackslash}m{2.5cm}}
Space & $k=0$, $d=2$ &  $k=1$, $d=2$ & $k=0$, $d=3$ & $k=1$, $d=3$ \\
\hline
& & & & \\
$ V_\rho$
& \includegraphics[scale=0.2]{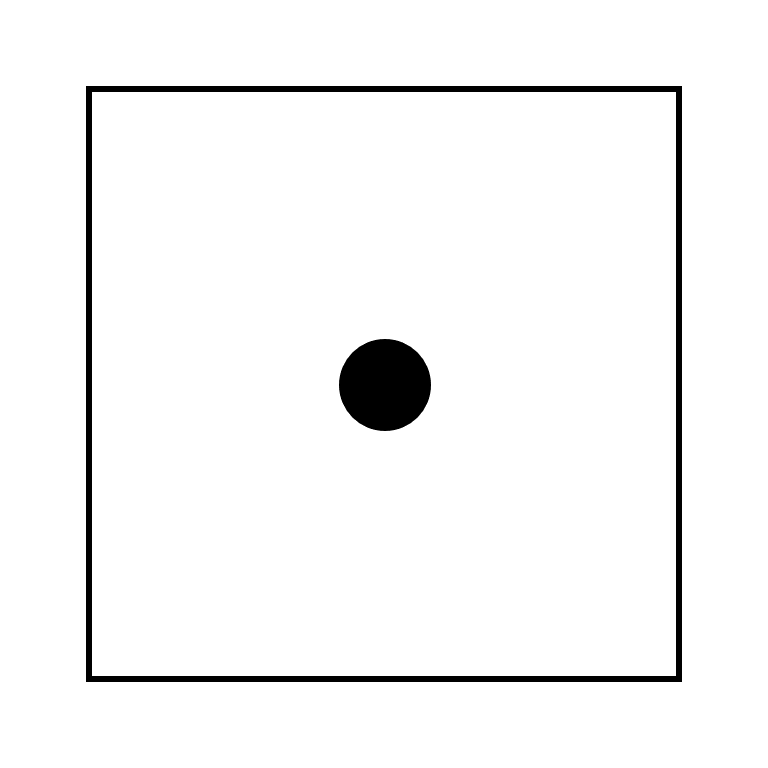}
& \includegraphics[scale=0.2]{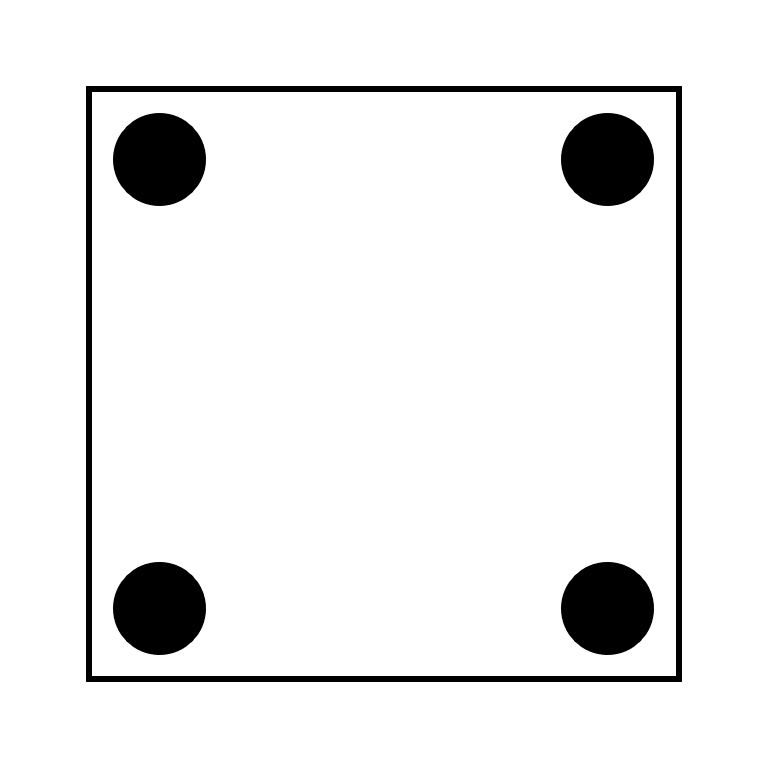}
& \includegraphics[scale=0.2]{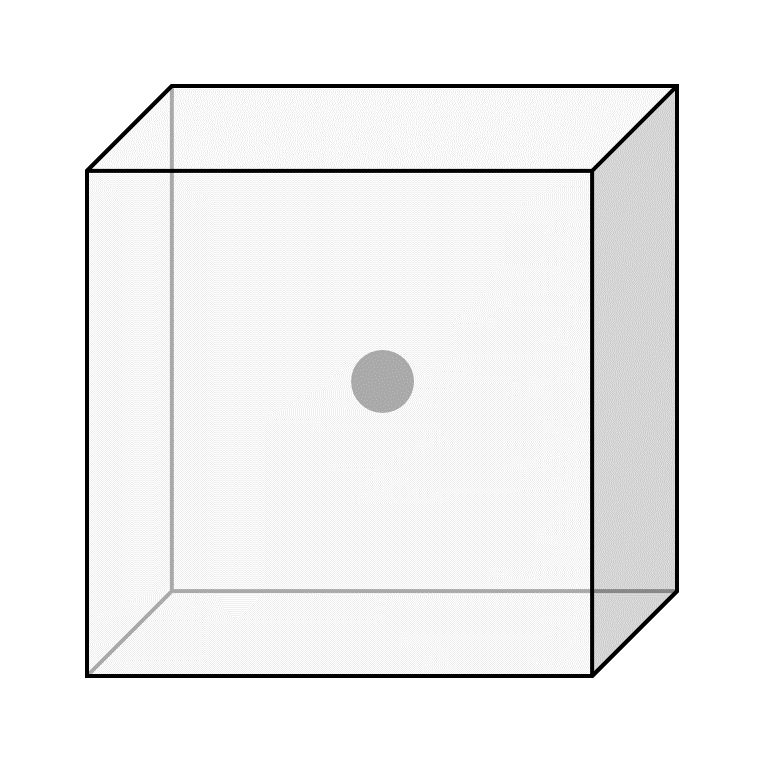}
& \includegraphics[scale=0.2]{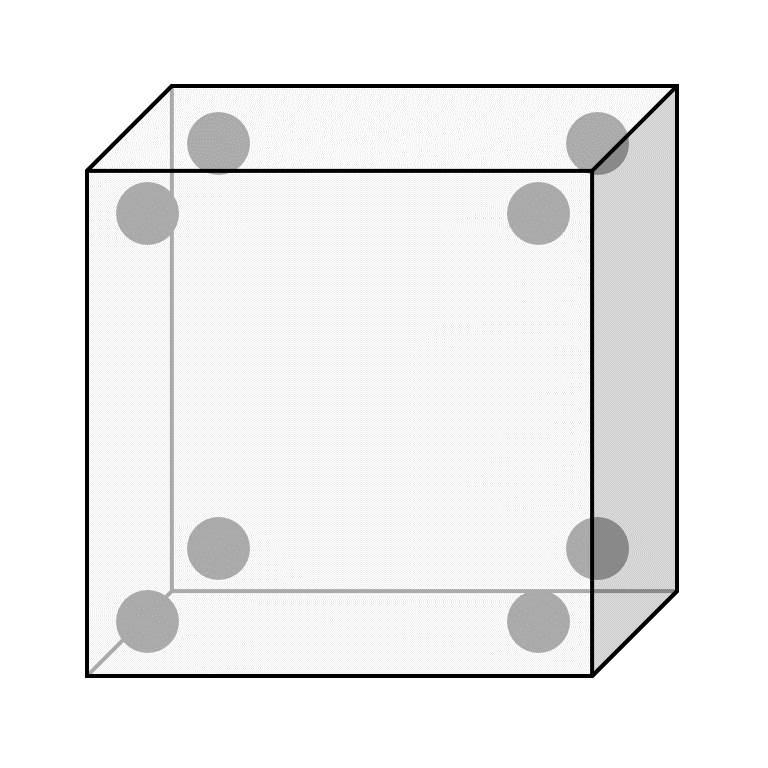}
\\
\hline
& & & & \\
$ V_{\bm{v}}$ 
& \includegraphics[scale=0.2]{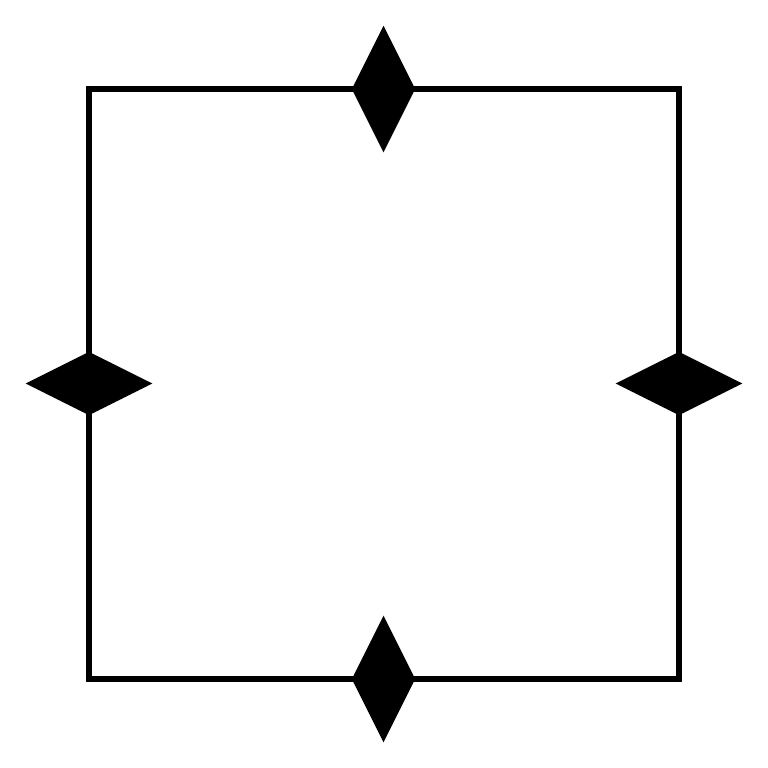}
& \includegraphics[scale=0.2]{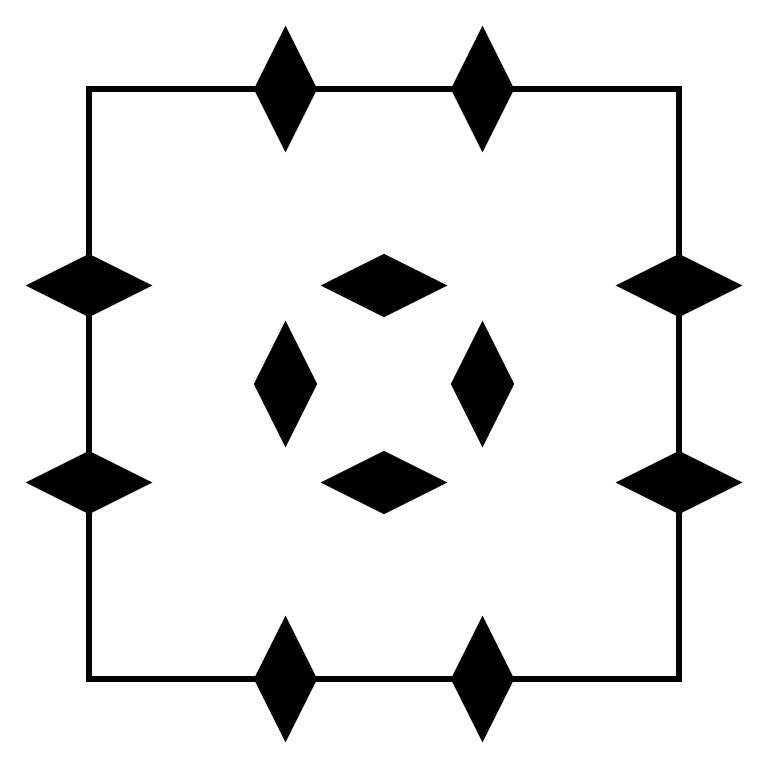}
& \includegraphics[scale=0.2]{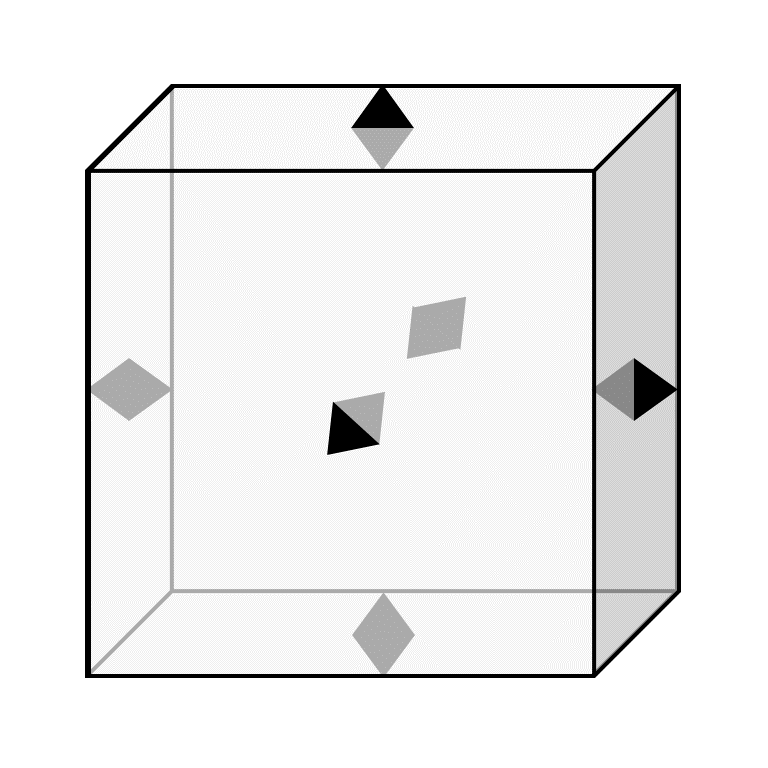}
& \includegraphics[scale=0.2]{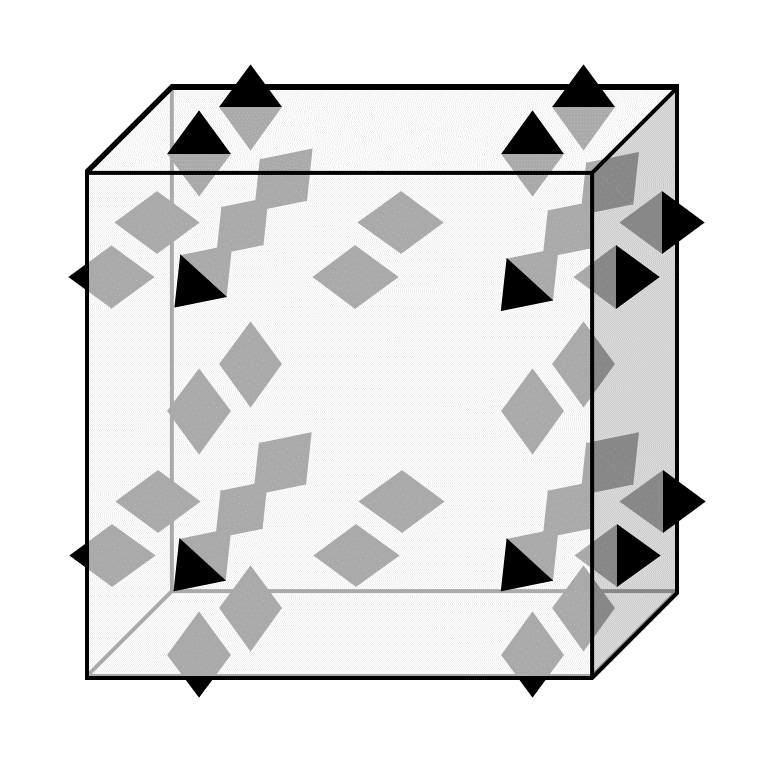}
\\
\hline
& & & & \\
$ V_\theta$ 
& \includegraphics[scale=0.2]{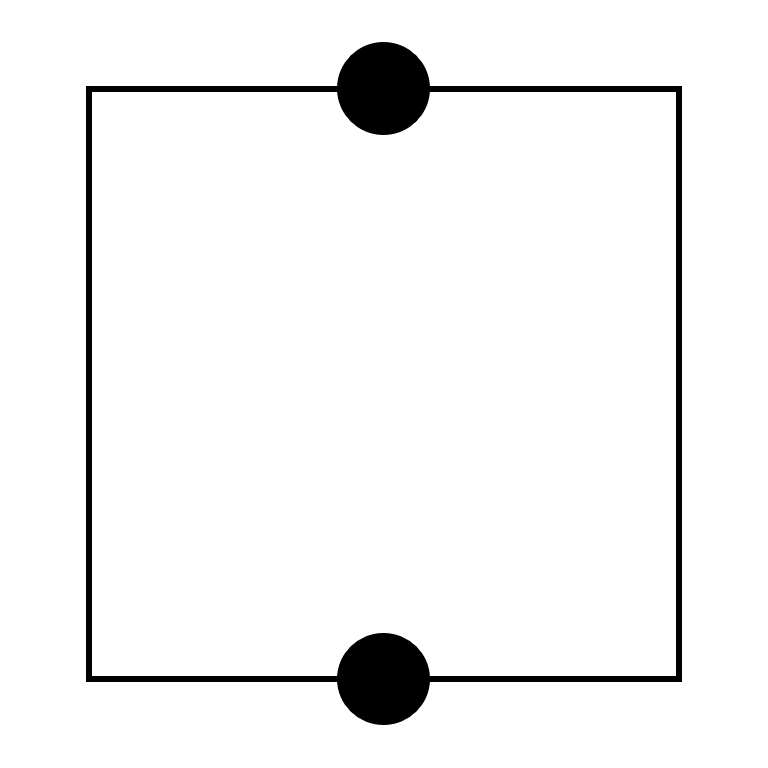}
& \includegraphics[scale=0.2]{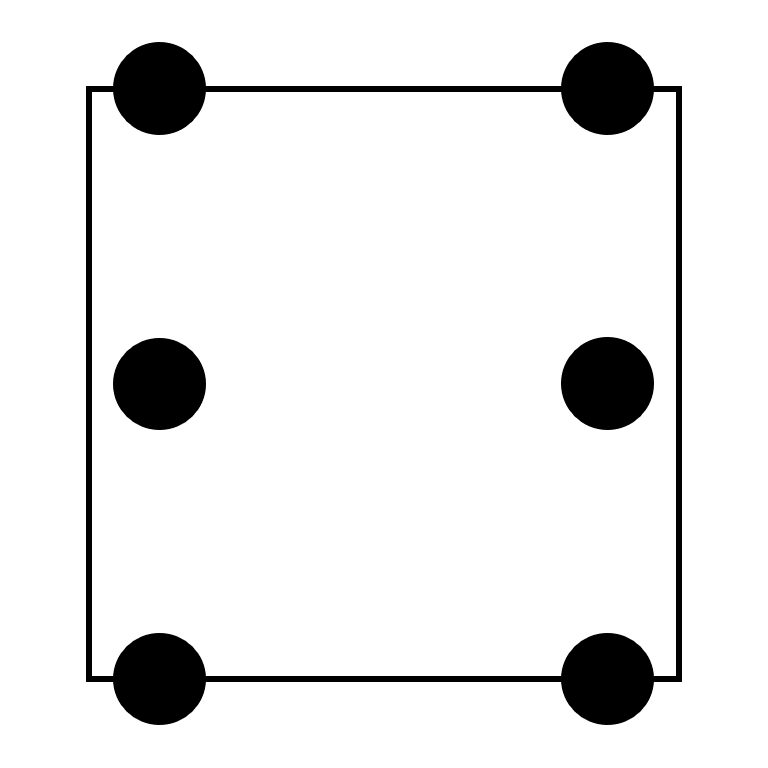}
& \includegraphics[scale=0.2]{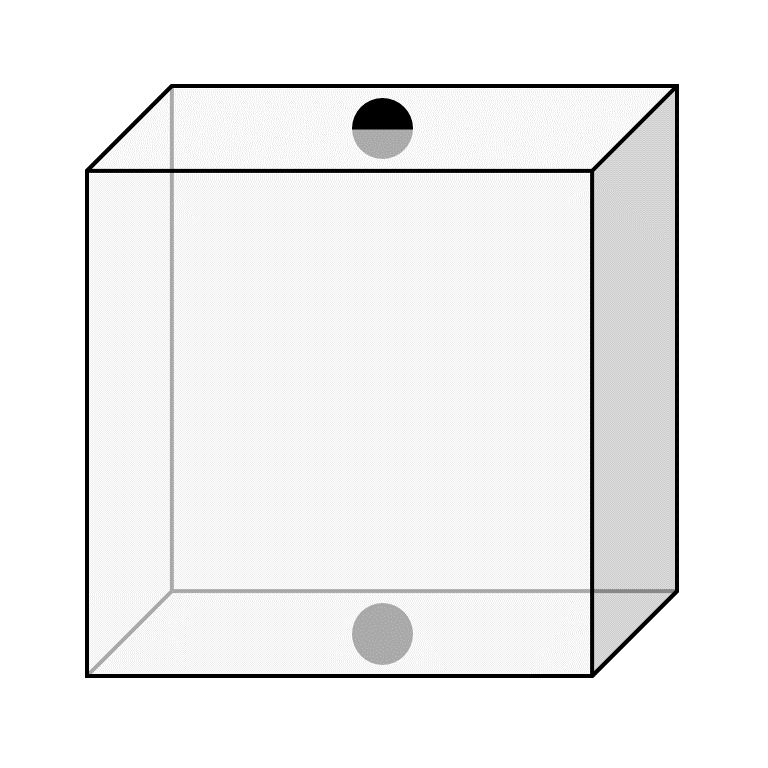}
& \includegraphics[scale=0.2]{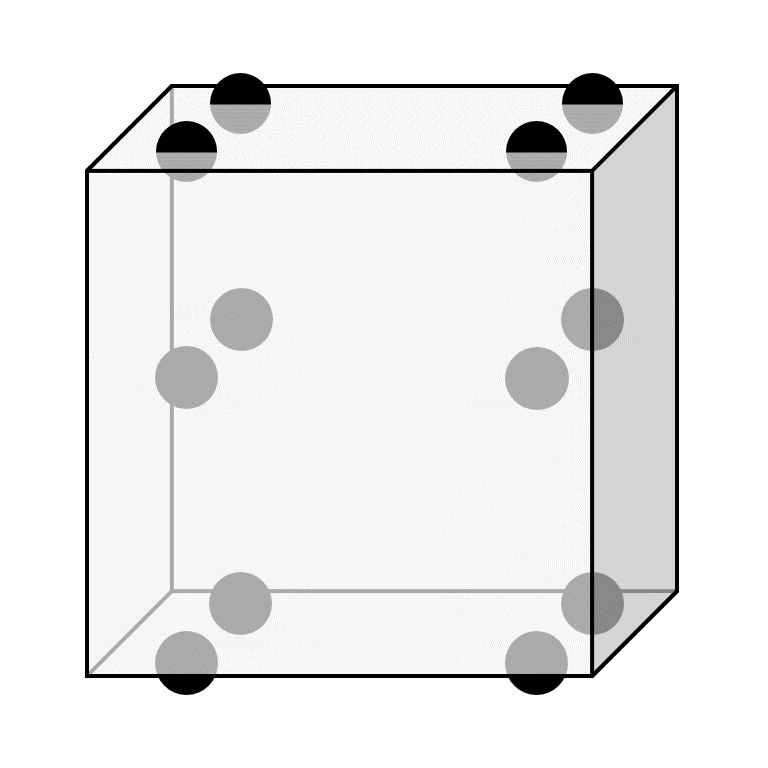}
\end{tabular}
\caption{An illustration of the function spaces that we use for density, velocity and potential temperature variables in vertical slice ($d=2$) and three-dimensional ($d=3$) simulations.
We have have two configurations for each case: with vertical and horizontal degree of $k=0$ and $k=1$.
Degrees of freedom on cell facets are shared between cells and denote the continuity of the field between them.} \label{tab:spaces}
\end{table}

\subsection{Overview}
Here we present an outline of our model structure. It is based the UK Met Office ENDGame model \cite{endgame2014documentation} in terms of the semi-implicit formulation and physics-dynamics coupling, with the exception that Eulerian advection schemes are used with Eulerian averages in the semi-implicit time discretisation. 
This is in contrast to the ENDGame semi-Lagrangian approach.
We illustrate this with some pseudocode that is similar to that presented in \cite{bendall2019recovered}.
The pseudocode describes how we evolve our prognostic variables, which we denote by the single vector $\bm{\chi}=(\bm{v}, \rho_d, \theta_{vd}, r_v, r_c, r_r)$.
Overall the model is a mixture of explicit and semi-implicit treatments.
Terms in \ref{eqn:wet atmos full} are divided into forcings, $\mathcal{F}(\bm{\chi})$, advection terms $\mathcal{A}(\bm{\chi})$ and physics terms $\mathcal{P}(\bm{\chi})$, with forcing being semi-implicit while advection and physics are explicit.\\
\\
Given $\bm{\chi}^n$, the state at the $n$-th time step, the first stage of the time step is to apply the explicit component of the forcing terms, returning $\bm{\chi^*}$.
Then we enter two loops.
In the outer loop, the advecting velocity $\overline{\bm{u}}$ is determined to be the average of $\bm{v}^n$ and the best guess of the velocity at the time step, $\bm{v}_p^{n+1}$.
The variables are then advected explicitly by $\overline{\bm{u}}$, returning $\bm{\chi}_p$.
To solve the implicit component of the forcing terms we enter the inner loop.
The residual $\Delta \bm{\chi}$ is calculated between $\bm{\chi}_p$ and the implicit part of the forcing, so that the implicit part of the model is solved if $\Delta \bm{\chi}=\bm{0}$.
We iterate towards the solution by solve a linearised system of equations, with the right-hand-side as $\Delta \bm{\chi}$.
Once the outer loop is completed, the physics processes are treated explicitly in their own stage.
\\
\\
The pseudocode summarising the algorithm is:
\begin{enumerate}[leftmargin=*]
\item \texttt{FORCING:}  $\bm{\chi}^*=\bm{\chi}^n+\tfrac{1}{2}\Delta t\mathcal{F}(\bm{\chi}^n)$
\item \texttt{SET:}   $\bm{\chi}_p^{n+1}=\bm{\chi}^n$
\item \texttt{OUTER:}
	\begin{enumerate}
	\renewcommand{\labelitemi}{}
	\item \texttt{UPDATE:}   $\overline{\bm{u}}=\tfrac{1}{2}(\bm{v}_p^{n+1}+\bm{v}^n)$
	\item \texttt{ADVECT:} $\bm{\chi}_p=\mathcal{A}_{\bar{\bm{u}}}(\bm{\chi}^*)$
	\item \texttt{INNER:}
	\begin{enumerate}
	\item \texttt{FIND RESIDUAL:} $\Delta\bm{\chi}= 
	\bm{\chi}_p+\tfrac{1}{2}\Delta t \mathcal{F}(\bm{\chi}_p^{n+1})
	-\bm{\chi}_p^{n+1}	$
	\item \texttt{SOLVE:} $\mathcal{S}(\bm{\chi}') = 
	\Delta\bm{\chi}$ \texttt{for} $\bm{\chi}'$
	\item \texttt{INCREMENT:} $\bm{\chi}_p^{n+1}=\bm{\chi}_p^{n+1}+\Delta\bm{\chi}$
	\end{enumerate}
	\end{enumerate}
\item \texttt{PHYSICS:} $\bm{\chi}_p^{n+1} = \mathcal{P}(\bm{\chi}_p^{n+1})$
\item \texttt{ADVANCE TIME STEP:} $\bm{\chi}^n = \bm{\chi}^{n+1}_p$
\end{enumerate}

\subsubsection{Differences between $k=0$ and $k=1$ set-ups}
Our model is designed for use with both the $k=0$ and $k=1$ sets of function spaces.
The formulation is slightly different between these two cases, and here we give a summary of those differences.
These are listed in Table \ref{tab:model configurations}.\\
\\
The main area of difference is in the advection part of the model, which is briefly detailed in Section \ref{sec:advection}.
In the $k=0$ case, we use the recovered space transport schemes outlined in \cite{bendall2019recovered}.
However when using the $k=1$ spaces, a discontinuous Galerkin upwind scheme is used for the transport of $\rho_d$, whilst the advection of $\theta_{vd}$ and the moisture variables is performed by the embedded scheme introduced by \cite{cotter2016embedded}.
The velocity $\bm{v}$ advection equation is written in vector invariant form, and uses the theta time stepping method.
All other advection schemes use a three-step Runge-Kutta time stepping procedure (SSPRK-3) \cite{shu1988efficient}. \\
\\
The other difference is that the limiters used with the advection schemes are different between the two set-ups.

\begin{table}
\centering
\begin{tabular}{|c|c|c|}
\hline & & \\ [\dimexpr-\normalbaselineskip+2pt]
\textbf{Component} &  $\bm{k=0}$ \textbf{configuration} &  $\bm{k=1}$ \textbf{configuration} \\
\hline
 Forcing & \multicolumn{2}{c|}{Non-advective terms of equations} \\
\hline
Transport of $\bm{v}$ &
\begin{tabular}{c}
Recovered DG \\
upwinding, SSPRK-3 \\
time stepping
\end{tabular} & 
\begin{tabular}{c} 
Vector invariant \\
form, theta method \\
time stepping
\end{tabular} \\
\hline & & \\ [\dimexpr-\normalbaselineskip+2pt]
Transport of $\rho_d$ & 
\begin{tabular}{c}
Recovered DG \\
upwinding, SSPRK-3  \\
time stepping
\end{tabular} &
\begin{tabular}{c}
DG upwinding, \\
SSPRK-3 \\
time stepping
\end{tabular} \\
\hline & & \\ [\dimexpr-\normalbaselineskip+2pt]
\begin{tabular}{c}
Transport of $\theta_{vd}$ \\
and moisture
\end{tabular} & 
\begin{tabular}{c}
Recovered DG \\
upwinding, SSPRK-3 \\
time stepping
\end{tabular} & 
\begin{tabular}{c}
Embedded DG \\ 
upwinding, SSPRK-3 \\
time stepping
\end{tabular} \\ 
\hline 
Linear solve & \multicolumn{2}{c|}{Hybridised linear solver} \\
\hline
\end{tabular}
\caption[Gusto model configurations]{A summary of the components of the different model configurations when using Gusto with the lowest order $k=0$ spaces, and the next-to-lowest order $k=1$ spaces.} \label{tab:model configurations}
\end{table}
\subsection{Forcing}\label{sec:forcing}
The `forcing' operation constitutes the sum of the non-advective terms in the differential equations (\ref{eqn:wet atmos momentum}) and (\ref{eqn:wet atmos theta_v eqn of motion}), which are written in weak form.
Our goal is to find $\mathcal{F}(\bm{v})$ and $\mathcal{F}(\theta_{vd})$. \\
\\
For both $k=0$ and $k=1$ configurations, the forcing that we apply to $\bm{v}$ is the solution $\bm{v}_\mathrm{trial}\in V_{\bm{v}}$, for all $\bm{\psi}\in V_{\bm{v}}$, to
\begin{equation}
\begin{split}
\int_\Omega \bm{\psi}\bm{\cdot}\bm{v}_\mathrm{trial} \dx{x}
 = & \int_\Omega \left[c_{pd}\Pi\bm{\nabla}\bm{\cdot}\left(\frac{\theta_{vd}\bm{\psi}}{1+r_t}\right)
- \bm{\psi}\bm{\cdot}(\bm{f}\times\bm{v})
\right] \dx{x} \\
& -\int_\Omega g \left(\bm{\psi\cdot}\widehat{\bm{k}}\right)\dx{x} - \int_\Gamma c_{pd} \left\llbracket \frac{\theta_{vd} \bm{\psi}}{1+r_t} \right\rrbracket_{\bm{n}} \langle \Pi \rangle \dx{S},
\end{split}
\end{equation}
where $\Omega$ is the domain, $\Gamma$ is the set of all interior facets and the angled brackets $\langle\cdot\rangle$ 
denote the average value on either side of a facet.
The double square brackets $\llbracket \bm{f} \rrbracket_{\bm{n}}=\bm{f}^+\bm{\cdot}\widehat{\bm{n}}^++\bm{f}^-\bm{\cdot}\widehat{\bm{n}}^-$ denote the jump in values between either side of the facet, arbitrarily labelled with $+$ and $-$.
The forcing term $\mathcal{F}(\bm{v})$ is then taken to the solution $\bm{v}_\mathrm{trial}$ to this. \\
\\
The forcing $\mathcal{F}(\theta_{vd})$ applied to $\theta_{vd}$ is the solution $\theta_\mathrm{trial}$, for all $\gamma\in V_\theta$ to
\begin{equation}
\int_\Omega \gamma \theta_\mathrm{trial}\dx{x}  = 
- \int_\Omega \gamma\theta_{vd}\left(\frac{R_m}{c_{vml}}-\frac{R_dc_{pml}}{c_{pd}c_{vml}}\right)(\dv{v})\dx{x}.
\end{equation}

\subsection{Advection}\label{sec:advection}
In the advection stage, each of the variables is translated by the velocity $\overline{\bm{u}}$.
We represent action this by the operator $\mathcal{A}_{\bar{\bm{u}}}$.
Here we briefly provide details of the advection schemes used.
For a summary, see Table \ref{tab:model configurations}.
\subsubsection{Discontinuous Galerkin Upwind Advection} \label{sec:upwind}
First, we define a single forward-Euler step of discontinuous Galerkin upwinding, which we describe as the operation $\mathcal{L}_{\bar{\bm{u}}}$ upon a scalar quantity $q$.
For the advective form of the transport equation this involves finding the solution $q_\mathrm{trial}$, for all $p\in V_q$
\begin{equation}
\int_\Omega \psi q_\mathrm{trial} \dx{x} -
\int_\Omega \psi q \dx{x} -
\Delta t\int_\Omega q \left[\grad{}\bm{\cdot} \left(p \overline{\bm{u}}\right)\right] \dx{x}
+ \Delta t\int_\Gamma \left(\bm{u}\bm{\cdot}\widehat{\bm{n}}^+ \right) q^* \llbracket p\rrbracket_{+} \dx{S}=0, \label{eqn:adv DG upwind}
\end{equation}
where the jump is $\llbracket p\rrbracket_{+} = p^+-p^-$.
The outward normal on the $+$ side of the facet is denoted by $\widehat{\bm{n}}^+$ and $q^*$ is the upwind value identified by
\begin{equation}
q^* = \left \lbrace \begin{matrix}
q^+, & \mathrm{if} \ \overline{\bm{u}}\bm{\cdot}\widehat{\bm{n}}^+ > 0, \\
q^-, & \mathrm{otherwise}.
\end{matrix}\right.
\end{equation} 
This equation can also be cast in continuity form, as it will be for the transport of $\rho_d$.
The SSPRK-3 time stepping scheme presented in \cite{shu1988efficient} involves composing these steps as follows:
\begin{subequations}\label{eqn:SSPRK3}
\begin{align}
q^{(1)} & := q+\mathcal{L}_{\bar{\bm{u}}}q , \\
q^{(2)} & := \frac{3}{4}q+\frac{1}{4}\left(q^{(1)}+\mathcal{L}_{\bar{\bm{u}}}q^{(1)}\right), \\
\mathcal{A}_{\bar{\bm{u}}}q & := \frac{1}{3}q + \frac{2}{3}\left(q^{(2)}+\mathcal{L}_{\bar{\bm{u}}}q^{(2)}\right).
\end{align}
\end{subequations}
This is the scheme used for transport for $\rho_d$ in the $k=1$ configuration, and is the basis for the embedded DG and recovered advection schemes.
\subsubsection{Embedded DG Advection}\label{sec:embedded DG}
The embedded DG transport scheme introduced in \cite{cotter2016embedded} is used for advecting $\theta_{vd}$ and the moisture variables in the $k=1$ configuration.
The fields are first injected into fully discontinuous analogues of their spaces, where the advection takes place.
The fields are then projected back to their original spaces.
See \cite{cotter2016embedded} for more details.

\subsubsection{Recovery Operator and Recovered Advection}\label{sec:recovery}
The advection schemes used for the $k=0$ spaces were presented in \cite{bendall2019recovered}.
This scheme extends the embedded DG scheme described in the previous section: the fields are `recovered' from the original space to spaces with higher-degrees of polynomial.
The transport then occurs in a discontinuous higher-degree space.
As in the embedded scheme, this field is projected back to the original space to complete the advection.
The spaces used as part of this scheme were described in \cite{bendall2019recovered}.\\
\\
The `recovery' operator is also used in physics-dynamics coupling, which is described in Section \ref{sec:parametrisations}.
The basic operation recovers a field from some function space $V_0$ to $\widetilde{V}_1$ by pointwise evaluation of the field at degrees of freedom.
When $V_0$ is discontinuous between elements but $\widetilde{V}_1$ is continuous, the new value of the field is the average over the neighbouring elements.
At the exterior boundaries of the domain, this process will not have second-order numerical accuracy and a further extrapolation step is required to accurately represent gradients.
\subsubsection{Vector invariant advection}
The most major difference in the advection stage of the model between using the $k=0$ and $k=1$ spaces is the transport of the velocity.
We write the advection equation in vector-invariant form, discretising
\begin{equation}
\pfrac{\bm{v}}{t} + \left(\grad{}\times \bm{v}\right) \times \overline{\bm{u}} + \frac{1}{2}\grad{}\left(\bm{v}\bm{\cdot}\overline{\bm{u}}\right)
=\bm{0}.
\end{equation}
The action of the advection operator $\mathcal{L}_{\bar{\bm{u}}}$ gives the solution $\bm{v}_\mathrm{trial}\in V_{\bm{v}}$ for all $\bm{\psi}\in V_{\bm{v}}$ to
\begin{equation}
\begin{split}
\int_\Omega \bm{\psi}\bm{\cdot}\left(\bm{v}_\mathrm{trial}-\bm{v}\right) \dx{x} 
& 
 -  \int_\Gamma \bm{v}^* \bm{\cdot}
\left[ \widehat{\bm{n}}^+\times \llbracket\overline{\bm{u}} \times \bm{\psi}\rrbracket_{+}\right]\mathrm{d}S \\
& + \int_\Omega \bm{v} \bm{\cdot} \left[\grad{}\times\left(\overline{\bm{u}}\times \bm{\psi}\right) \right] \mathrm{d}x 
- \frac{1}{2}\int_\Omega \left(\bm{v}\bm{\cdot}\overline{\bm{u}}\right)\left(\grad{}\bm{\cdot}
\bm{\psi}\right)\mathrm{d}x=0,
\end{split}
\end{equation}
with $\bm{v}^*$ as the upwind value of $\bm{v}$ and $\llbracket \bm{\cdot}\rrbracket_{+}$ taking the same definition as in Section \ref{sec:upwind}.  \\
\\
The theta time stepping scheme is used, so that the velocity at the $(n+1)$-th time step is found by solving
\begin{equation}
\bm{v}_{n+1} = \bm{v}_n + \Delta t \left[\vartheta \mathcal{L}_{\bar{\bm{u}}}\bm{v}_n + (1-\vartheta)\mathcal{L}_{\bar{\bm{u}}}\bm{v}_{n+1} \right]
\end{equation}
We take $\vartheta = 1/2$.

\subsubsection{Limiting}
\begin{figure}
\centering
{\includegraphics[width=0.48\textwidth]{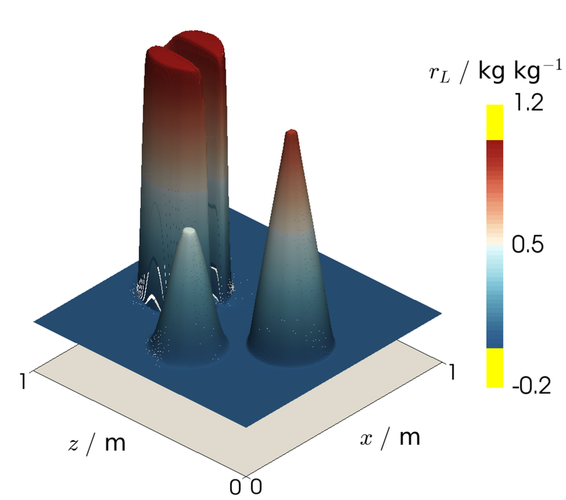}}
{\includegraphics[width=0.48\textwidth]{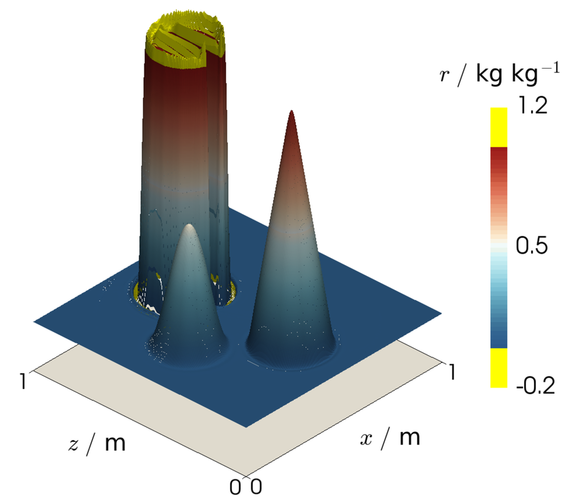}}
\caption[Effectiveness of limiter in the $k=1$ configuration]{The final fields, following one rotation around the domain, of functions in the two-dimensional $\mathrm{dQ}_1\otimes\mathrm{Q}_2$ space that were initialised with the slotted-cylinder, hump, cone initial condition of \cite{leveque1996high}.
(Left) with the limiter applied, and (right) without the limiter. Overshoots and undershoots are highlighted in yellow.} \label{fig:k1 limiter}
\end{figure}

The advection schemes that we use do not preserve the monotonicity property of the continuous transport equation.
This can result in unphysical negative concentrations when advecting moisture variables.
A way of preventing this is to apply a slope limiter. \\
\\
We provide the option of adding a limiter to the advection of variables in $V_\theta$.
The vertex-based limiter of \cite{kuzmin2010vertex} is designed for use on fields that are piecewise-linear and discontinuous between cells.
It prevents the formation of new maxima and minima by separating the field into a constant mean part and adjusting the perturbed part in each cell.
It is applied to the field before the advection operator $\mathcal{V}_{\bar{\bm{u}}}$ first acts each time step and after each stage of the three-step Runge Kutta method. \\
\\
For advection of moisture variables in the $k=0$ set of spaces, we use the recovered space scheme and the advection can be limited by the vertex-based limiter of \cite{kuzmin2010vertex}.
When projecting back from the advecting space $V_1$, as described in \cite{bendall2019recovered} we split the projection into two steps: projecting into the broken space $\widehat{V}_0$ before recovering any continuity. \\
\\
With the $k=1$ spaces, the vertex-based limiter of \cite{kuzmin2010vertex} is again used.
However this field is quadratic in the vertical and the values of the field at vertical midpoint of the cell are unchanged by the vertex-based limiter.
We bound these values by restricting them to the average of the values at the adjacent vertices if it falls outside of the range spanned by them.
In other words, if $r_{i+1/2}$ is the value of the field at a degree of freedom halfway up the $i$-th cell, then the new value is
\begin{equation}
r_{i+1/2} = \left\lbrace \begin{matrix}
r_{i+1/2} & \mathrm{if} \ \min\left(r_{i}, r_{i+1}\right) \leq r_{i+1/2} \leq \max\left(r_{i}, r_{i+1}\right), \\
\frac{1}{2}\left(r_{i}+ r_{i+1} \right), & \mathrm{otherwise}.
\end{matrix}\right.
\end{equation}
To illustrate the effectiveness of this limiter, we performed the rotation of the slotted-cylinder, hump, cone of \cite{leveque1996high} in a 2D vertical slice.
Figure \ref{fig:k1 limiter} shows the fields in the $V_\theta$ space for the $k=1$ configuration, rotated once around the domain by a velocity defined by the stream function
\begin{equation}
\psi = \pi(z(z-1)+x(x-1)).
\end{equation}
The domain used was the unit square, with $\Delta x=\Delta z=0.01$, whilst the time step was $\Delta t=10^{-4}$.
The final fields were recorded at $t=1$.
\subsection{Solve} \label{sec:solver}
We now present the strategy for the linear solve step of the model.
When initialising the model, we set mean $\rho_d$ and $\theta_{vd}$ states which we denote respectively as $\overline{\rho}_d$ and $\overline{\theta}_{vd}$.
The residual $\Delta\bm{\chi}$ between the predicted explicit $\bm{\chi}_p$ and the latest guess $\bm{\chi}_p^{n+1}$ of the prognostic variables is computed first.
Reducing this residual to zero solves the implicit part of the model.
This is done by solving a linear set of equations, with the residual on the right hand side, to calculate increments that are added to $\bm{\chi}_p^{n+1}$, updating it.\\
\\
The linearised problem that we are attempting to solve is
\begin{subequations}
\begin{align}
\bm{v}'+\frac{\Delta t c_{pd}}{2(1+r_t)}\left(\theta'_{vd}\grad{\overline{\Pi}}+\overline{\theta}_{vd}\grad{\Pi}'\right) & = \Delta \bm{v}, \\
\rho_d'+\tfrac{1}{2}\Delta t\grad{}\bm{\cdot}\left(\overline{\rho}_d\bm{v}'\right) & = \Delta \rho_d, \\
\theta_{vd}'+\tfrac{1}{2} \Delta t \left(\widehat{\bm{k}}\bm{\cdot}\bm{v}'\right)\left(\widehat{\bm{k}}\bm{\cdot}\grad{\overline{\theta}_{vd}}\right) & = \Delta \theta_{vd},
\end{align}
\end{subequations}
where the primes represent the perturbations to be found. \\
\\
These equations will be solved in weak form.
Howver, we simplify them first by eliminating $\theta_{vd}'$.
This allows us to create a hybridised mixed system to be solved for $\bm{v}'_\mathrm{trial}$ and $\rho'_\mathrm{trial}$.
Taking $\widehat{\bm{k}}$ as the unit upward normal, we introduce
\begin{equation}
\Theta' = \Delta\theta_{vd} - \tfrac{1}{2} \Delta t \left(\widehat{\bm{k}}\bm{\cdot}\grad{\overline{\theta}_{vd}}\right)\left(\widehat{\bm{k}}\bm{\cdot}\bm{v}'_{\mathrm{trial}}\right).
\end{equation}
The Exner pressure perturbation is then approximated as
\begin{equation}
\Pi' = \frac{\kappa}{1-\kappa} \overline{\Pi}\left(\frac{\Theta'}{\overline{\theta}_{vd}} + \frac{\rho'_\mathrm{trial}}{\overline{\rho}_d}\right), \ \ \mathrm{where} \  \
\overline{\Pi}=\left(\frac{\overline{\rho}_d \overline{\theta}_{vd} R_d}{p_R}\right)^\frac{\kappa}{1-\kappa}.
\end{equation}
The hybridised mixed system involves breaking the continuity of normal components of $\bm{v}'$ between elements, introducing the broken space $\widehat{V}_{\bm{v}}$.
In conjunction, a trace field $\ell'_\mathrm{trial}$ is introduced approximating the average values of $\Pi$ perturbations on the trace space $V_\mathrm{trace}$: the discontinuous space of functions that live on the boundary of $V_{\bm{v}}$.
This variable acts as Lagrange multipliers to provide the continuity constraint of the pre-hybridised system.
The resulting mixed system for $(\bm{v}'_\mathrm{trial},\rho'_\mathrm{trial},\ell'_\mathrm{trial})\in \left(\widehat{V}_{\bm{v}}, V_\rho, V_\mathrm{trace}\right)$ which holds for all
$\bm{\psi}\in \widehat{V}_{\bm{v}}$, $\phi\in V_\rho$ and $\lambda\in V_\mathrm{trace}$, can be written:
\begin{equation} \label{eqn:moist solve}
\begin{split}
& \int_\Omega \bm{\psi}\bm{\cdot}\left(\bm{v}'_\mathrm{trial}
-\Delta\bm{v}\right) \dx{x} \\
 & - \frac{\Delta t}{2} c_{pd}\left(\int_\Omega \grad{}\bm{\cdot}\left[\frac{\Theta' w}{1+r_t}\widehat{\bm{k}}\right] \overline{\Pi}\dx{x} 
-  \int_\Gamma \left\llbracket 
\frac{\Theta' w}{1+r_t}\widehat{\bm{k}} \right\rrbracket_{\bm{n}} \left\langle \overline{\Pi} \right\rangle \dx{S}\right) \\
& + \frac{\Delta t}{2} c_{pd}\left(\int_{\partial \Omega} 
\frac{\Theta' w}{1+r_t}\widehat{\bm{k}}\bm{\cdot}\widehat{\bm{n}}  \left\langle \overline{\Pi} \right\rangle \dx{s}
-\int_\Omega \grad{}\bm{\cdot}\left[\frac{\overline{\theta}_{vd}}{1+r_t}\bm{\psi}\right]\Pi'\dx{x}\right) \\
& + \frac{\Delta t}{2} c_{pd}\left(\int_{\Gamma} 
\left\llbracket\frac{\overline{\theta}_{vd}}{1+r_t}\bm{\psi} \right\rrbracket_{\bm{n}} \ell'_\mathrm{trial} \dx{S}
+\int_{\partial\Omega} \frac{\overline{\theta}_{vd}}{1+r_t}\left(\bm{\psi}\bm{\cdot}\widehat{\bm{n}}\right)\ell'_\mathrm{trial}\dx{s}\right) \\
& +\int_\Omega \phi \left(\rho'_\mathrm{trial}-\Delta\rho_d\right)\dx{x}  - \frac{\Delta t}{2} \int_\Omega \left(\grad{\phi}\bm{\cdot}\bm{v}'_\mathrm{trial}\right)\overline{\rho}_d \dx{x} \\
& + \frac{\Delta t}{2} \left(\int_\Gamma \left\llbracket \phi \bm{v}'_\mathrm{trial} \right\rrbracket_{\bm{n}} \left\langle \rho_d\right \rangle\dx{S}
+ \int_{\partial\Omega}\phi \bm{v}'_\mathrm{trial}\bm{\cdot}\widehat{\bm{n}} \left\langle \rho_d \right \rangle \dx{s} \right) \\
& + \int_\Gamma \lambda \left\llbracket \bm{v}'_\mathrm{trial}\right\rrbracket_{\bm{n}} \dx{S}
+ \int_{\partial \Omega}\lambda \left(\bm{v}'_\mathrm{trial}\bm{\cdot}\widehat{\bm{n}}\right)\dx{s} =0.
\end{split}
\end{equation}
where $w=\widehat{\bm{k}}\bm{\cdot}\bm{\psi}$ is the vertical component of the velocity test function and $\partial \Omega$ is the external boundary of the domain.
This mixed system can be statically condensed into a single system for the Lagrange multipliers, giving a single system to be solved.
Once this system has been solved, the calculated value of $\ell'$ is used to find $\rho'_d$ and then the broken $\bm{v}'$.
The recovery operator (averaging values between cells) is then used to restore the continuous velocity $\bm{v}'$.
Finally, the value of $\theta_{vd}$ is found as the $\theta'_\mathrm{trial}$ that solves for all $\gamma\in V_\theta$
\begin{equation}
\int_\Omega \gamma \left[\theta'_\mathrm{trial} - \Delta\theta_{vd} +\frac{\Delta t}{2} \left(\widehat{\bm{k}}\bm{\cdot}\bm{v}\right)
\left(\widehat{\bm{k}}\bm{\cdot}\grad{}\overline{\theta}_{vd}\right] \right)=0.
\end{equation}
More on this solver and its performance will be detailed in future work.

\section{Physics Parametrisations} \label{sec:parametrisations}
In this section we consider the discretisation of the `physics' terms in equation set (\ref{eqn:wet atmos full}). 
Each term is evaluated individually in a separate `physics' routine.
The term labelled $\dot{r}^c_\mathrm{cond}$ represents the condensation of water vapour into cloud water and the evaporation of cloud water into water vapour.
We use the rate given in \cite{rutledge1983mesoscale} and \cite{bryan2002benchmark}, with the saturation mixing ratio stemming from Tetens' empirical formula \cite{tetens1930formula}, giving
\begin{equation}
r_\mathrm{sat}(\widetilde{p},T)=\cfrac{\epsilon C^\mathrm{sat}_0}{\widetilde{p}\exp\left[-\cfrac{C^\mathrm{sat}_1(T-T_R)}{T-C^\mathrm{sat}_2}\right]-C^\mathrm{sat}_0}, \label{eqn:rsat}
\end{equation}
where $C^\mathrm{sat}_0$, $C^\mathrm{sat}_1$ and 
$C^\mathrm{sat}_2$ are constants whose values are given in the appendix.
The pressure $\widetilde{p}$ will be explained in the following section. \\
\\
The similar term $\dot{r}^r_\mathrm{evap}$ describes the  evaporation of rain water into water vapour, with the value taken from \cite{klemp1978simulation}. \\
\\
The coalescence of cloud water into rain droplets is described by $\dot{r}_\mathrm{accr}$ and $\dot{r}_\mathrm{accu}$: the accretion and auto-accumulation processes respectively.
The formulas that we use are also those described in \cite{klemp1978simulation}, with values from \cite{soong1973comparison}.\\
\\
The sedimentation of rain is given by $S$.
Our approach to determining $S$ is similar
to the single moment scheme described in \cite{milbrandt2010sedimentation}.
We assume that the number $n_r(D)$ of raindrops of diameter $D$ forms a spectrum described by a Gamma distribution.
This is used to determine the terminal velocity of raindrops averaged over $D$.
The rain field $r_r$ is then advected by this averaged velocity using the same advection scheme used by the dynamics for $r_r$ described in Section \ref{sec:dynamics}.
\subsection{Combining Fields From Different Function Spaces}
Many diagnostic fields, such as temperature $T$ or saturation mixing ratio of water vapour $r_\mathrm{sat}$ must be determined from prognostic fields in different function spaces, usually $\rho_d$, $\theta_{vd}$ and $r_v$.
Such diagnostic fields are important in the equations used for the physics parametrisations and can also be important in setting the initial state of the model. \\
\\
As $\rho_d$ and $\theta_{vd}$ lie in different function spaces, when determining some diagnostic variable $q$ there are a number of different choices that can be taken.
With our emphasis on parametrisations of moist processes, we take only the case of $q\in V_\theta$, so that it is also in the same function space as the moisture variables.
Our approach to combining these fields is to `recover' $\rho_d$ is into $V_\theta$ using the recovery operator outlined in Section \ref{sec:recovery} to give $\widetilde{\rho}_d$.
Then $q$ can be calculated algebraically within $V_\theta$.
This approach can be used with both the $k=0$ and the $k=1$ sets of spaces, and with the $k=0$ spaces it has second-order numerical accuracy (including at the boundaries of the domain).
For example, the temperature $T$ and pressure $p$ used in parametrisations are calculated as
\begin{equation}
T = \frac{\theta_{vd}\widetilde{\Pi}}{1+r_v/\epsilon}, \ \ \ \ \
\widetilde{p}=p_R\widetilde{\Pi}^\frac{1}{\kappa},
\end{equation}
where 
\begin{equation}\label{eq:exner-pi}
\widetilde{\Pi} = \left(\frac{\widetilde{\rho_d}R_d\theta_{vd}}{p_R}
\right)^\frac{\kappa}{1-\kappa}.
\end{equation}
\subsection{Time Discretisation}
At present we are performing simple explicit first-order
integration in time for the physics routines.
In other words for a process $\dot{r}$ affecting a variable $r$,
the new value $r_\mathrm{new}$ will be related to the old 
$r_\mathrm{old}$ by
\begin{equation}
r_\mathrm{new}=r_\mathrm{old}+\Delta t \ \dot{r}(r_\mathrm{old}).
\end{equation}
The exception here is the treatment of the sedimentation of
rainfall, which was described in the previous section.
The value of $\dot{r}$ comes from the state of the model
just before this physics routine is called.
This is the state of the model after either the dynamics
or the preceding physics routine has been completed.
The physics routines are executed consecutively, (this is commonly known as ``sequential splitting''), in the following order:
\begin{enumerate}
\item accretion of cloud water,
\item auto-accumulation of rain water,
\item sedimentation of rain,
\item evaporation of rain water,
\item evaporation of cloud water/condensation of water vapour.
\end{enumerate}
We choose to do the evaporation/condensation step last so as to
prevent any supersaturation at the end of the time step.
\section{Hydrostatic Balance Routines}\label{sec:hydrostatic balance}
For many test cases the background or initial state of the model will be in hydrostatic balance.
The procedure that we use to obtain a discrete hydrostatic balance is based upon that presented in \cite{natale2016compatible}.
Given a boundary condition for the pressure and a $\theta_d$ field, this procedure finds the $\rho_d$ which gives rise to zero vertical velocity. 
In this section we describe two developments to this: the extension of the routine of \cite{natale2016compatible} to use a hybridised formulation and the treatment of cases in which thermodynamic and moist variables and the prognostic variables need to be found.
\subsection{Hybridised Hydrostatic Balance}
The hybridised hydrostatic balance system uses a similar approach to that described in Section \ref{sec:solver}.
As in that case, the variables are expressed in discontinuous spaces and Lagrange multipliers $\ell$ are introduced in a trace space on the horizontal facets to provide continuity constraints.
We also introduce $\mathring{V}_{\bm{v}}$, the subspace of the broken velocity space $\widehat{V}_{\bm{v}}$ with zero flow normal to the boundaries.
The problem is then to find the $(\bm{v}_\mathrm{trial},\rho_\mathrm{trial},\ell_\mathrm{trial})\in \left(\mathring{V}_{\bm{v}}, V_\rho, V_\mathrm{trace}\right)$, such that for all $(\bm{\psi},\phi,\lambda)\in \left(\mathring{V}_{\bm{v}}, V_\rho, V_\mathrm{trace}\right)$:
\begin{subequations}
\begin{align} \label{eqn:hybrid-hydrostatic}
\begin{split}
\int_\Omega \bm{\psi\cdot} \bm{v}_\mathrm{trial}\dx{x} -  \int_{\Omega} 
c_{pd}\widehat{\Pi}(\rho_d,\theta_{vd}) & \grad{}\bm{\cdot}\left(\frac{\theta_{vd}\bm{\psi}}{1+r_t}\right)\dx{x} \\
+ \int_\Gamma \llbracket \bm{\psi}\rrbracket_{\bm{n}}\dx{S} & = - \int_{\partial\Omega_0}\frac{c_{pd}\theta_{vd}}{1+r_t}\left(\bm{\psi}\bm{\cdot}\widehat{\bm{n}}\right)\Pi_0\dx{S} - \int_\Omega g \bm{\psi}\bm{\cdot}\widehat{\bm{k}}\dx{x},
\end{split} \\
\int_\Omega \phi c_{pd} \grad{}\bm{\cdot}\left(\frac{\theta_{vd}\bm{v}_\mathrm{trial}}{1+r_t}\right)\dx{x} & = 0, \\
\int_\Gamma \lambda \llbracket \bm{v}_\mathrm{trial}\rrbracket_{\bm{n}} \dx{S} & = 0, 
\end{align}
\end{subequations}
where $\partial\Omega_0$ is the boundary upon which the condition that $\Pi=\Pi_0$ is to be satisfied and $\widehat{\Pi}$ is a linearisation of $\Pi$. 
Equation (\ref{eqn:hybrid-hydrostatic}) can be statically condensed to an elliptic equation for $\ell$, which can be efficiently inverted in each column of the mesh.
In fact it turns out that $\ell$ is an approximation for $c_{pd}\theta_{vd}\Pi/(1+r_t)$ on the horizontal facets of the mesh.
Once $\ell$ is determined, $\bm{v}$ and $\rho_d$ can be solved for via local back-substitution in each cell of the mesh.
Continuity is restored in $\bm{v}$ again by using the recovery operator.
By hybridising, we avoid solving a larger, more ill-conditioned mixed system and instead inverted a condensed elliptic problem.
The implementation of this procedure is described by \cite{gibson2019slate}, with the application taking the form of a custom Python-base preconditioner conforming to standard PETSc library options \cite{petsc1997,petty2011modified}, as described in \cite{gibson2019slate}.
The advantage being that solver options for solving the condensed
system arising from reducing \eqref{eqn:hybrid-hydrostatic} can easily be updated, even when
nested inside a non-linear method.

\subsection{Saturated Conditions}\label{sec:sat balance}
This set-up involves initial conditions like those in the moist
benchmark of \cite{bryan2002benchmark}.
The problem is to find $\rho_d$, $\theta_{vd}$ and $r_v$ given
$\theta_e$, $r_t$ and a boundary condition on the pressure.
Assuming the absence of rain, $r_v=r_\mathrm{sat}$,
and $r_c = r_t - r_v$.
The \textit{wet-equivalent potential temperature}, $\theta_e$,
is a conserved quantity in reversible, moist adiabatic processes, 
i.e. $\mathrm{D}\theta_e/\mathrm{D}t=0$.
Following \cite{emanuel1994atmospheric}, for our equation set 
(\ref{eqn:wet atmos full}), $\theta_e$ can be written as
\begin{equation}\label{eqn: theta_e def}
\theta_e := T\left(\frac{p_0}{p}\right)^\frac{R_d}{c_{pd} + c_{pl}r_t}
\left(\mathcal{H}\right)^\frac{-r_vR_v}{c_{pd}+c_{pl}r_t}
\exp\left[\frac{L_v(T) r_v}{(c_{pd} + c_{pl}r_t)T}\right], 
\end{equation}
which with $\mathcal{H}=1$ for a saturated atmosphere 
is the same as that used by 
\cite{bryan2002benchmark} and that derived in the appendix of 
\cite{paluch1979entrainment} from the second law of
thermodynamics.
The saturation mixing ratio is given by (\ref{eqn:rsat}). \\
\\
The challenge is to obtain the $\theta_{vd}$ and $\rho_d$ fields
that satisfy the specified $\theta_e$ field whilst ensuring
that $r_v = r_\mathrm{sat}$.
We use an initial guess for $\theta_{vd}$ and feed it to equation
(\ref{eqn:hybrid-hydrostatic}), to generate a
guess for $\rho_d$.
This density is then converted into $V_\theta$, before a nested fixed-point iteration-style procedure is
used to invert $\theta_e(\theta_{vd}, \widetilde{\rho}_d, r_v)$
and $r_\mathrm{sat}(\theta_{vd}, \widetilde{\rho}_d, r_v)$,
to obtain $\theta_{vd}$ and $r_v$. \\
\\
Let $l$, $m$ and $n$ count the number of iterations to find $\rho_d$, $\theta_{vd}$ and $r_v$ respectively.
These form nested loops, such that the outer loop uses the latest approximations of $\theta_{vd}$ and $r_v$ with equation (\ref{eqn:hybrid-hydrostatic}) to obtain $\rho_d^{(l+1/2)}$.
We then combine this with the previous value to update the approximation to $\rho_d$, doing
\begin{equation}
\rho_d^{(l+1)}=(1-\delta)\rho_d^{(l)} + \delta \rho_d^{(l+1/2)},
\end{equation}
where $\delta=0.8$.
The next loop finds $\theta_{vd}$ using
\begin{equation}
\theta^{(m+1/2)}_{vd} = \frac{\theta^{(m)}_{vd} \theta_e}{\theta_e\left(\theta^{(m)}_{vd}, \widetilde{\rho}^{(l)}_d, r_v^{(n)}\right)}, \ \ \ 
\theta_{vd}^{(m+1)}=(1-\delta)\theta_{vd}^{(m)} + \delta \theta_{vd}^{(m+1/2)}
\end{equation}
whilst the inner loop computes
\begin{equation}
r_v^{(n+1/2)} = r_\mathrm{sat}\left(\theta^{(m)}_{vd}, \widetilde{\rho}^{(l)}_d, r_v^{(n)}\right), \ \ \
r_v^{(n+1)}=(1-\delta)r_v^{(m)} + \delta r_v^{(m+1/2)}
\end{equation}
where $\theta_e$ without arguments denotes the specified value.
These loops are iterated until
$\theta_e(\theta_{vd}, \widetilde{\rho}_d, r_v)$ converges to the
specified value to some tolerance.
\subsection{Unsaturated Conditions} \label{sec:unsat balance}
We now discuss how to find the prognostic thermodynamic variables given $\theta_d$ and the relative humidity $\mathcal{H}$, such as the case in Section \ref{sec:unsat bubble} or \cite{grabowski1991cloud}.
The relative humidity is related to $r_v$ by
\begin{equation}
\mathcal{H}=\frac{r_v}{r_\mathrm{sat}}\left(\frac{1+r_\mathrm{sat}/ \epsilon}{1+r_v / \epsilon}\right). \label{eqn:RH}
\end{equation}
As in Section \ref{sec:sat balance}, we use what can be thought of as a nested iterative procedure.
Counting the latest approximations of $\rho_d$, $\theta_{vd}$ and $r_v$ with $l$, $m$ and $n$, the outer loop uses $\theta_{vd}^{(m)}$ and $r_v^{(n)}$ in the hydrostatic balance equation (\ref{eqn:hybrid-hydrostatic}) to determine $\rho_d^{(l+1/2)}$.
Again, the next value of $\rho_d$ is given by
\begin{equation}
\rho_d^{(l+1)}=(1-\delta)\rho_d^{(l)} + \delta \rho_d^{(l+1/2)}.
\end{equation}
There is only one inner loop in this case, so that $m=n$.
The new value of $r_v$ is found by rearranging \ref{eqn:RH} so that
\begin{equation}
r_v^{(m+1/2)} = \frac{\mathcal{H} r_\mathrm{sat}\left(\theta^{(m)}_{vd}, \widetilde{\rho}^{(l)}_d, r_v^{(m)}\right)}{1 + (1-\mathcal{H})r_\mathrm{sat}\left(\theta^{(m)}_{vd}, \widetilde{\rho}^{(l)}_d, r_v^{(m)}\right) / \epsilon}, \ \ \ r_v^{(m+1)}=(1-\delta)r_v^{(m)} + \delta r_v^{(m+1/2)},
\end{equation}
where $\mathcal{H}$ is the specified value of the relative humidity.
The final step is to use the specified value of $\theta_d$ to get
\begin{equation}
\theta_{vd}^{(m)} =  \theta_d \left(1 + r_v^{(m)}/\epsilon\right).
\end{equation}
This iterative process continues until $\mathcal{H}\left(\theta_{vd}^{(m)},\rho_d^{(l)}, r_v^{(m)}\right)$ has converged to its specified value to some tolerance.

\section{Test Cases} \label{sec:test cases}
In this section we demonstrate the discretisation detailed in previous sections through a series of test cases, with some comparison of the $k=0$ and $k=1$ configurations of the model.
Two new test cases are presented, featuring a gravity wave in a saturated atmosphere and a three-dimensional rising thermal in a saturated atmosphere.
\subsection{Bryan and Fritsch Moist Benchmark}\label{sec:sat bubble}
\begin{figure}[h!]
\centering
{\includegraphics[width=0.4\textwidth]{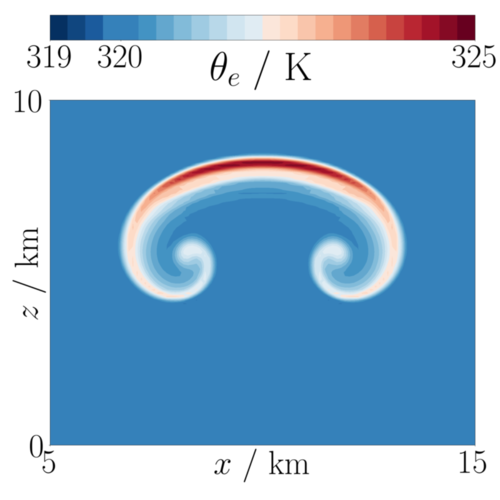}} ~~
{\includegraphics[width=0.4\textwidth]{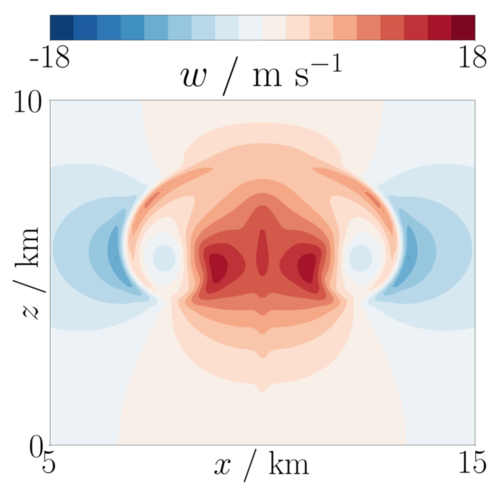}}
\caption[Two dimensional saturated thermal with $k=0$ spaces]{ The (left) $\theta_e$ field contoured every $0.25$ K and (right) vertical velocity $w$ field contoured every $2$ m s$^{-1}$, with both fields plotted at $t=1000$ s for a simulation of the moist benchmark case \cite{bryan2002benchmark} representing a thermal rising through a saturated atmosphere.
The $320$ K contour has been omitted for clarity in the $\theta_e$ field.
This simulation is with the $k=0$ lowest-order set of spaces, with grid spacing $\Delta x=\Delta z=100$ m and a time step of $\Delta t=1$ s.
These solutions are visibly similar to those presented in \cite{bryan2002benchmark}.} \label{fig:moist bf 1}
\end{figure}

\begin{figure}[h!]
\centering
{\includegraphics[width=0.4\textwidth]{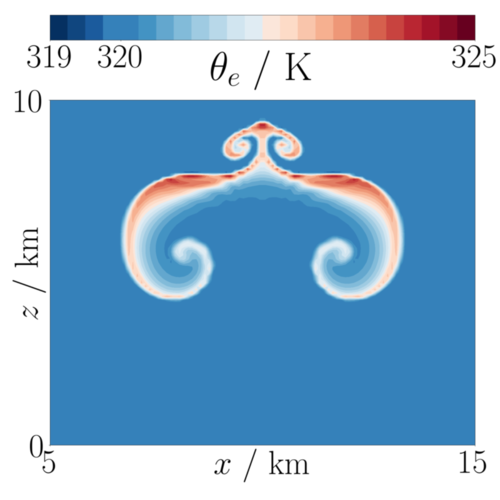}} ~~
{\includegraphics[width=0.4\textwidth]{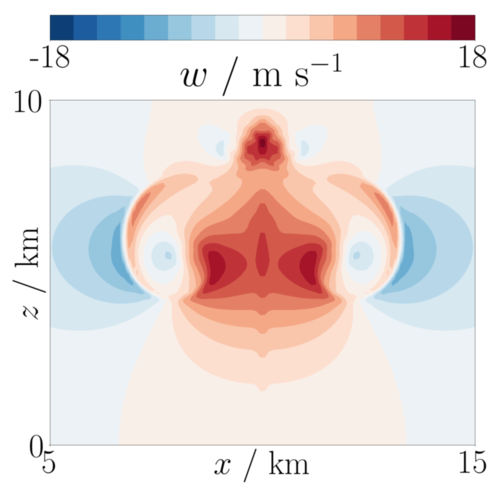}}
\caption[Two dimensional saturated thermal with $k=1$ spaces]{Outputted fields from the $k=1$ next-to-lowest degree space simulation at $t=1000$ s of the moist benchmark from \cite{bryan2002benchmark}. (Left) $\theta_e$ with contours spaced by $0.25$ K and (right) vertical velocity $w$ contoured every $2$ m s$^{-1}$.
The simulation used grid spacing $\Delta x=\Delta z=100$ m and a time step of $\Delta t=1$ s.
The $320$ K contour has been omitted for clarity in the $\theta_e$ field.
A second plume can be seen forming at the top of the primary plume.}\label{fig:moist bf 2}
\end{figure}

The first demonstration of our discretisation is through the moist benchmark test case of \cite{bryan2002benchmark}, which simulates a rising thermal through a cloudy atmosphere.
The domain is a vertical slice of width $L=20$ km and height $H=10$ km.
Periodic boundary conditions are applied at the vertical boundaries, but the top and bottom boundaries are rigid, so that $\bm{v\cdot}\widehat{\bm{n}}=0$ along them.
As in \cite{bryan2002benchmark}, we include no rain microphysics and no Coriolis force.\\
\\
The initial conditions defined in \cite{bryan2002benchmark} specify a background state with constant $r_t=0.02$ kg kg$^{-1}$ and constant wet-equivalent potential temperature $\theta_e=320$ K, which is defined in (\ref{eqn: theta_e def}).
Along with these, the background state is given by the requirements of hydrostatic balance, $r_v=r_\mathrm{sat}$ and $p=10^5$ Pa at the bottom boundary.
The procedure described in Section \ref{sec:sat balance} allows us to obtain the prognostic variables $\theta_{vd}$, $\rho_d$, $r_v$ and $r_c$ that approximately satisfy these conditions. \\
\\
A perturbation is then applied to $\theta_{vd}$.
With $(x,z)$ as the horizontal and vertical coordinates, the perturbed field is
\begin{equation}
\theta'_{vd} = \left\lbrace 
\begin{matrix}
\Delta\Theta \cos^2\left(\frac{\pi r}{2r_c}\right), & r<r_c, \\
0, & \mathrm{otherwise},
\end{matrix} \right.
\end{equation}
where $\Delta \Theta=2$ K, $r_c=2$ km and with $x_c=L/2$ and $z_c=2$ km we define
\begin{equation}
r = \sqrt{\left(x-x_c\right)^2 + \left(z-z_c\right)^2}.
\end{equation}
The initial $\theta_{vd}$ field is given in terms of the background field $\overline{\theta}_{vd}$:
\begin{equation}
\theta_{vd} = \overline{\theta}_{vd} \left(1 + \frac{\theta_{vd}'}{300 \ \mathrm{K}} \right).
\end{equation}
In the test case of \cite{bryan2002benchmark}, the pressure field is unchanged by the perturbation.
To replicate this with our prognostic variables, we find $\rho_d$ such that for all $\zeta\in V_\rho$,
\begin{equation}
\int_\Omega \zeta \rho_d \theta_{vd} \dx{x} = \int_\Omega \zeta \overline{\rho}_d\overline{\theta}_{vd} \dx{x},
\end{equation}
where $\overline{\rho}_d$ and $\overline{\theta}_{vd}$ are the background states for $\rho_d$ and $\theta_{vd}$.
The system is returned to saturation by finding the $r_v$ that solves, for all $\phi \in V_{\theta}$,
\begin{equation}
\int_\Omega \phi r_v \dx{x} = \int_\Omega \phi  r_\mathrm{sat}(\widetilde{\rho}_d, \theta_{vd}, r_v) \dx{x},
\end{equation}
where $r_\mathrm{sat}(\widetilde{\rho}_d, \theta_{vd}, r_v)$ is the expression for saturation mixing ratio in terms of the initial $\rho_d$ recovered into $V_\theta$ and the initial $\theta_{vd}$ field that have already been found, and also the unknown $r_v$ to be solved for.
Finally, $r_c$ is found from applying $r_c=r_t-r_v$.
The initial velocity field is zero in each component. \\
\\
Figures \ref{fig:moist bf 1} and \ref{fig:moist bf 2} show the $\theta_e$ and vertical velocity $w$ fields at $t=1000$ s.
Figure \ref{fig:moist bf 1} shows these fields for the configuration using the $k=0$ lowest-order spaces, with the $k=1$ spaces shown in Figure \ref{fig:moist bf 2}.
Both simulations used $\Delta x=\Delta z = 100$ m and $\Delta t=1$ s.
These final states are visibly different: whilst the $k=0$ solutions resemble those of \cite{bryan2002benchmark}, the $k=1$ solution displays an extra plume forming at the top of the rising thermal.
We believe this to be a manifestation of a physical instability that is damped by numerical diffusion in the $k=0$ case.
The $k=1$ solution appears highly sensitive to the choice of mesh, as at higher resolution the $\theta_e$ field does not appear to converge to a single solution.
Indeed if the domain is spanned horizontally by an odd number of cells, rather than a secondary plume emerging, the top of the primary plume appears to collapse.
This behaviour is also observed in the absence of moisture.

\subsection{Inertia-Gravity Waves in Saturated Atmosphere}\label{sec:moist gravity waves}
\begin{figure}
\centering
{\includegraphics[width=0.4\textwidth]{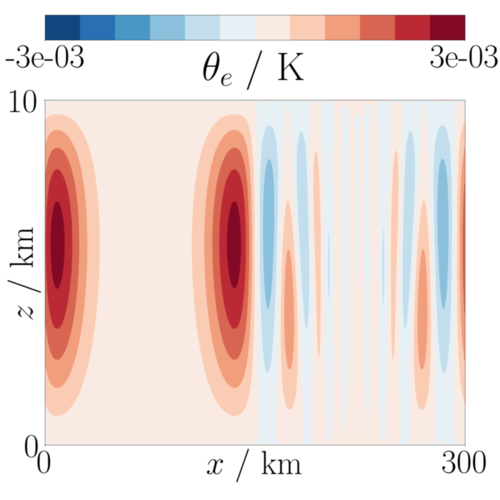}} ~~
{\includegraphics[width=0.4\textwidth]{moist_sk_LO_theta_e_pert.png}}
\caption[Gravity waves in a saturated atmosphere]{The perturbations to the $\theta_e$ fields at $t=3600$ s for the moist gravity wave test case. (Left) the $k=0$ lowest-order spaces set-up using $\Delta x=500$ m and (right) the $k=1$ spaces with $\Delta x=1000$ m. Both cases used $\Delta t=1.2$ s. Contours are spaced every $5\times 10^{-4}$ K.} \label{fig:moist gravity waves}
\end{figure}

We present here a new test case, a moist version of the non-hydrostatic gravity wave test of \cite{skamarock1994efficiency}, but in a saturated atmosphere like that of the moist benchmark in \cite{bryan2002benchmark}.
The final state of this test is spatially smooth, making this test appropriate for convergence tests.
No rain physics is used in this test case, and there is also no Coriolis force. \\
\\
The problem is set up in a two-dimensional vertical slice of length $L=300$ km and height $H=10$ km.
The dry gravity wave setup used in \cite{skamarock1994efficiency} applies a perturbation to $\theta_d$, which has a stratified background profile and is in hydrostatic balance.
Our variation on this is to apply the perturbation to the stratified background profile of $\theta_e$ in hydrostatic balance.
Using $(x,z)$ as the horizontal and vertical coordinates, the specified $\theta_e$ profile is
\begin{equation}
\overline{\theta}_e = \Theta_0 e^{N^2z/g},
\end{equation}
where $\Theta_0=$ 300 K and $N^2=10^{-4}$ s$^{-2}$.
With $r_t=0.02$ kg kg$^{-1}$ everywhere and the boundary condition of $p=10^5$ Pa at $z=0$, we use the hydrostatic balance procedure laid out in Section \ref{sec:sat balance} to find the $\theta_{vd}$, $\rho_d$, $r_v$ and $r_c$ fields that correspond to these initial conditions with the requirements of hydrostatic balance and that $r_v=r_\mathrm{sat}$ everywhere.
The initial velocity applied is $\bm{v}=(U,0)$ with $U=20$ m s$^{-1}$ describing uniform flow in the $x$-direction.
This defines each of the mean fields.\\
\\
A perturbation is then added, which is specified as
\begin{equation}
\theta'_e = \frac{\Delta \Theta}{1 + a^{-2}(x - L/2)^2} \sin\left(\frac{\pi z}{H}\right),
\end{equation}
with $a=5\times10^3$ m and $\Delta\Theta=0.01$ K.
The perturbed initial condition is then given by $\theta_e=\overline{\theta}_e+\theta'_e$.
Setting the new requirements that both $r_t$ and the pressure are unchanged by the addition of the perturbation, and that we still have $r_v=r_\mathrm{sat}$, defines the problem necessary to solve to find the initial prognostic fields.
We do this via a nested iterative process related to that described in Section \ref{sec:sat balance}.
In the outer loop we find $\rho_d^h$ such that for all $\zeta\in V_\rho$
\begin{equation}
\int_\Omega \zeta \rho_d^h \theta_{vd}^n \dx{x} = 
\int_\Omega \zeta \overline{\rho}_d \overline{\theta}_{vd} \dx{x},
\end{equation}
which is combined with the previous best estimate of $\rho^n_d$ to give $\rho^{n+1}_d=(1-\delta)\rho_d^n + \delta \rho^h_d$,
where $\delta=0.8$.
Nested inside this process are more damped iterations to find $\theta_{vd}$ and $r_v$, exactly as in Section \ref{sec:sat balance}. \\
\\
Figure \ref{fig:moist gravity waves} shows the perturbation to the final diagnostic $\theta_e$ field at $t=3600$ s for simulations with the $k=0$ lowest-degree spaces with $\Delta x=\Delta z=500$ m and the $k=1$ set of spaces, where $\Delta x=\Delta z=1000$ m, both with $\Delta t=1.2$ s.
These different cases are not visibly different from one another, and closely resemble the final state of the dry case from \cite{skamarock1994efficiency}.
\subsubsection{Convergence}
\begin{figure}[h!]
\centering
\includegraphics[width=0.9\textwidth]{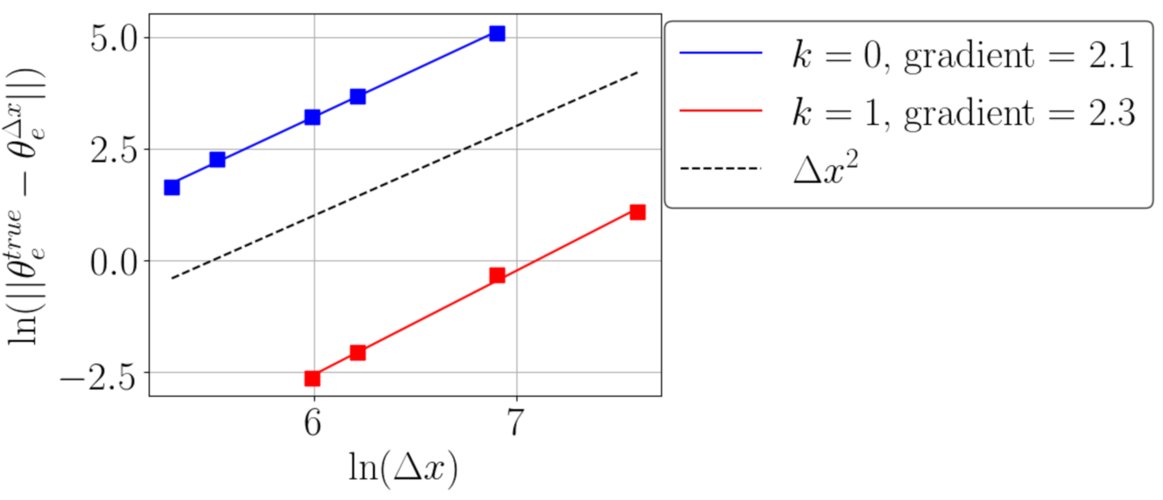}
\caption[Convergence of Gusto in solving moist equations]{A convergence plot showing the error as a function of resolution in the final state from the moist gravity wave test of Section \ref{sec:moist gravity waves}.
The true solution was taken from a high resolution simulation.
Both the $k=0$ and $k=1$ configurations of the model appear to have second-order or better accuracy.}
\label{fig:moist convergence}
\end{figure}

As the solution and evolution of this example is spatially smooth, we can use it to form a convergence test upon our model.
Although the final state (at $t=3600$ s) does not have an analytic solution, we use the fields from a high resolution simulation as the true solution. \\
\\
To measure the spatial accuracy of our model, we ran this test case at different resolutions, with each using a time step of $\Delta t=1.2$ s, which is small enough to avoid the Courant number breaching its critical value in the highest resolution cases.
The error is measured looking at the $\theta_e$ diagnostic in $V_\theta$ at $t=3600$ s.
The $\theta_e$ fields are interpolated onto the finest mesh, which has $\Delta x=100$ m for the $k=0$ case but $\Delta x=200$ m for the $k=1$ case.
The error between the high resolution solution for $\theta_e$ and those run at coarser resolutions is measured in the $L^2$ norm.
Results for our model are plotted in Figure \ref{fig:moist convergence}, which indicates that in both the $k=0$ and $k=1$ cases the model has second-order spatial accuracy, with the error proportional to $(\Delta x)^2$.

\subsection{Rising Thermal with Rain}\label{sec:unsat bubble}
This test case is based upon one described in \cite{grabowski1991cloud}.
This involves a thermal rising in an unsaturated atmosphere, forming a cloud which rains out.
This is another two dimensional vertical slice test, this time with domain of height $H=2.4$ km and length $L=3.6$ km, again with periodic conditions at the vertical sides.
The Coriolis force is neglected. \\
\\
In contrast to the saturated atmosphere initial conditions of Sections \ref{sec:sat bubble}, the initial state is defined by the dry potential temperature $\theta_d$ and a relative humidity field $\mathcal{H}$.
The background fields are $\mathcal{H}=0.2$ everywhere and
\begin{equation}
\theta_d = \Theta e^{Sz},
\end{equation}
where $\Theta$ is the dry potential temperature corresponding to $T_\mathrm{surf}=283$ K and $p=8.5\times10^4$ Pa, which also provides the pressure condition at the boundary.
The stratification is given by $S=1.3\times 10^{-5}$ m$^{-1}$.
We then use the procedure outlined in Section \ref{sec:unsat balance} to find the background $\theta_{vd}$, $\rho_d$ and $r_v$ fields that satisfy hydrostatic balance.
The initial $r_c$ and $r_r$ fields are zero. \\
\\
The perturbation is then applied to the relative humidity field $\mathcal{H}$, with a circular bubble that is just saturated, with an outer disk smoothing the perturbation into the background state.
This initial relative humidity field is given by
\begin{equation}
\mathcal{H} = \left \lbrace \begin{matrix}
\overline{\mathcal{H}}, & r \geq r_1, \\
\overline{\mathcal{H}} + (1 - \overline{\mathcal{H})}\cos^2\left(\frac{\pi(r-r_2)}{2(r_1-r_2)}\right), & r_2 \leq r < r_1, \\
1, & r<r_2,
\end{matrix} \right.
\end{equation}
where $\overline{\mathcal{H}}=0.2$ and
\begin{equation}
r = \sqrt{(x-x_c)^2 + (z-z_c)^2},
\end{equation}
with $x_c =L/2$, $z_c = 800$ m, $r_1=300$ m and $r_2=200$ m.
The $r_v$ and $\theta_{vd}$ that correspond to this $\mathcal{H}$ are found via a fixed point iterative method. \\
\begin{figure}[h!]
\centering
{\includegraphics[width=0.4\textwidth]{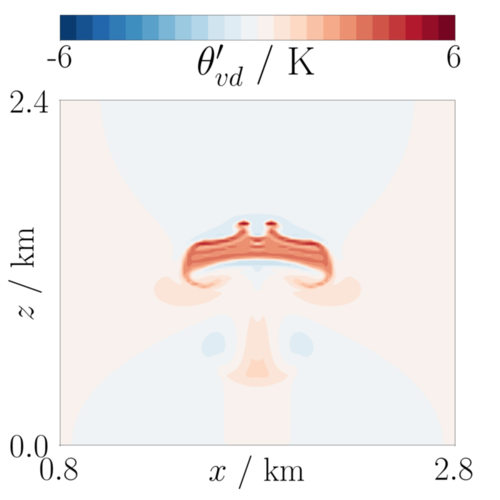}}~~
{\includegraphics[width=0.4\textwidth]{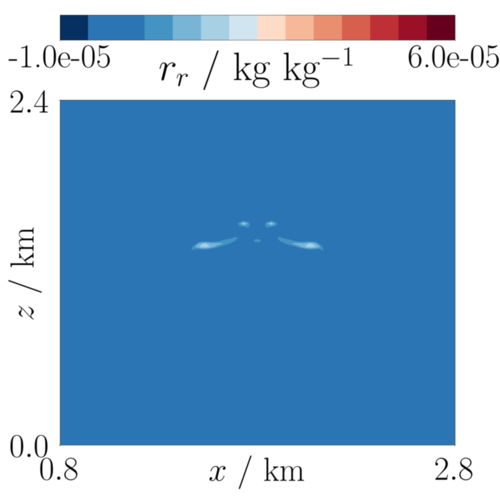}}
{\includegraphics[width=0.4\textwidth]{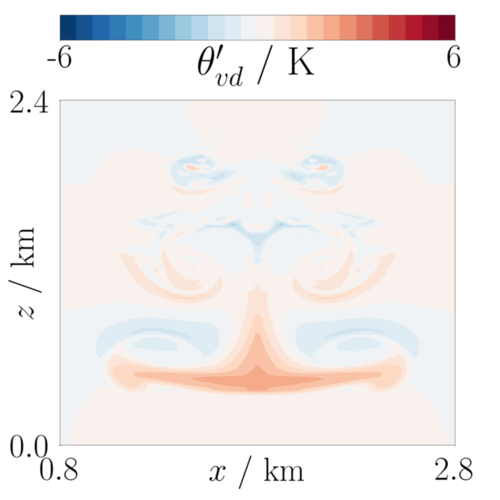}}~~
{\includegraphics[width=0.4\textwidth]{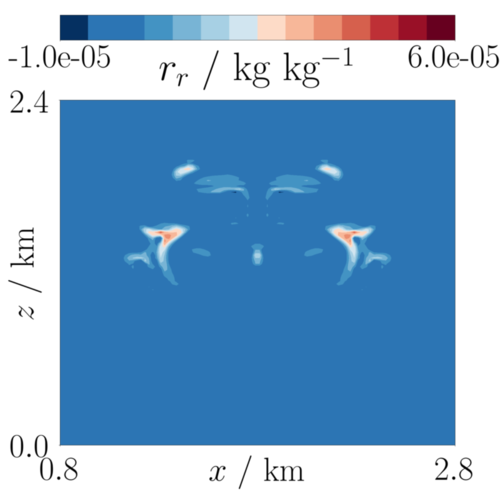}}
\caption[Rain bubble field output without limiter]{The field output at (top row) $t=300$ s and (bottom row) $t=600$ s from the rising thermal test with rain.
Fields shown are in the $k=0$ spaces with $\Delta x=\Delta z=20$ m and a time step of $\Delta t=1$ s.
(Left) the perturbation to $\theta_{vd}$, contoured every $0.25$ K.
(Right) the rain field $r_r$, with contours every $5\times 10^{-6}$ kg kg$^{-1}$ and omitting the zero contour.} \label{fig:rain results 1}
\end{figure}
\begin{figure}[h!]
\centering
{\includegraphics[width=0.4\textwidth]{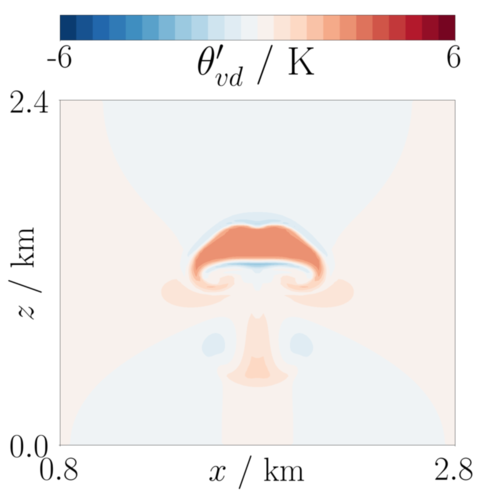}}~~
{\includegraphics[width=0.4\textwidth]{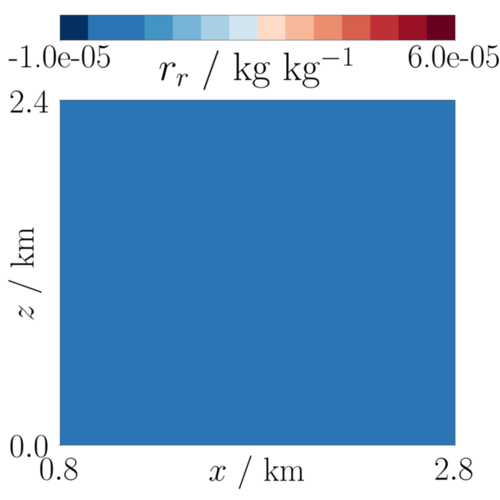}}
{\includegraphics[width=0.4\textwidth]{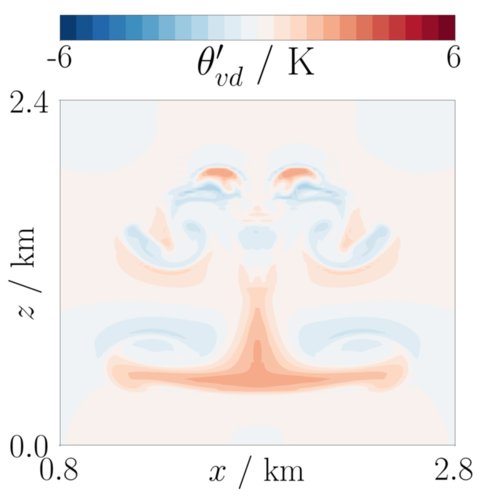}}~~
{\includegraphics[width=0.4\textwidth]{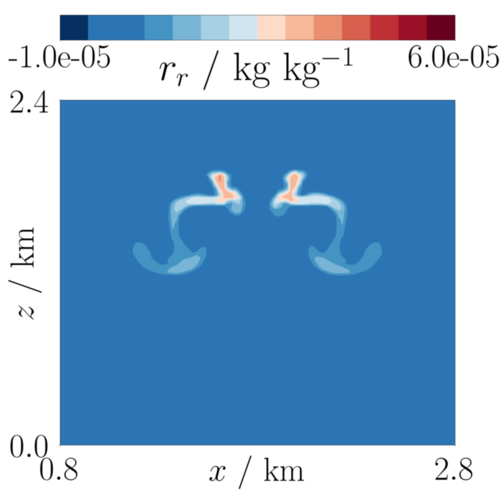}}
\caption[Rain bubble field output with limiter]{The same results as in Figure \ref{fig:rain results 1},  but with a limiter applied to the moisture species.} \label{fig:rain results 2}
\end{figure} \\
\noindent Fields are displayed in Figures \ref{fig:rain results 1} and \ref{fig:rain results 2} for $\theta'_{vd}$ and $r_r$ at $t=300$ and $600$ s.
Both simulations use the $k=0$ lowest-order space set-up, with $\Delta x=20$ m and $\Delta t=1$ s.
Figure \ref{fig:rain results 1} shows the results with no limiter applied to the advected moisture fields, whilst Figure \ref{fig:rain results 2} shows the same set-up but with a limiter applied to all the moisture variables.\\
\\
The length scales of the simulation are small enough for it to be highly turbulent, with the final state dependent on the resolution in the absence of a turbulence parametrisation.
Indeed, the lack of turbulence parametrisation in our model explains why these results look significantly different to those of \cite{grabowski1991cloud}.
Comparing Figures \ref{fig:rain results 1} and \ref{fig:rain results 2} demonstrates the effect of limiting the transport of moisture species.
In the absence of the limiter, the mass of water vapour depreciates less and so more cloud is formed, associated with a greater release of latent heat and a stronger updraught.
We see that rain forms earlier in the absence of a limiter.
However, some negative moisture values do form, which are absent from the limited case.
\subsection{Three-Dimensional Thermal in a Saturated Atmosphere}\label{sec:3d bubble}
\begin{figure}
\centering
{\includegraphics[width=0.4\textwidth]{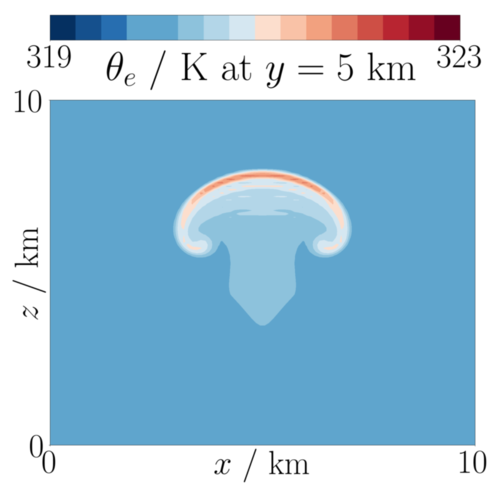}}~~
{\includegraphics[width=0.4\textwidth]{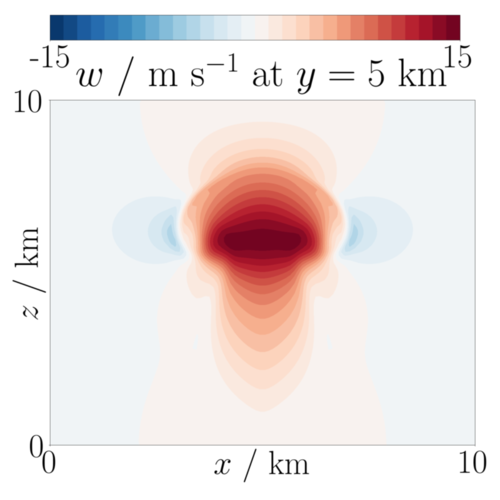}}
\caption[3D saturated rising bubble in $k=0$ configuration]{The (left) $\theta_e$ and (right) vertical velocity $w$ fields at $t=1000$ s for the three-dimensional simulation of a thermal rising through a saturated atmosphere.
The $\theta_e$ field is contoured every $0.25$ K with the $320$ K contour omitted, whilst the contour spacing for the $w$ field is $1$ m s$^{-1}$.
Cross-sections are shown at $y=5$ km, with values plotted on the lower side of the plane.
This simulation is with the $k=0$ lowest-order set of spaces, with grid spacing $\Delta x=100$ m and a time step of $\Delta t=1$ s.} \label{fig:3d bubble 1}
\end{figure}

\begin{figure}
\centering
{\includegraphics[width=0.4\textwidth]{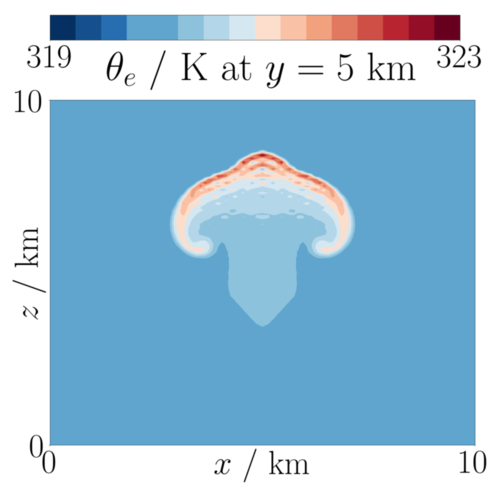}}~~
{\includegraphics[width=0.4\textwidth]{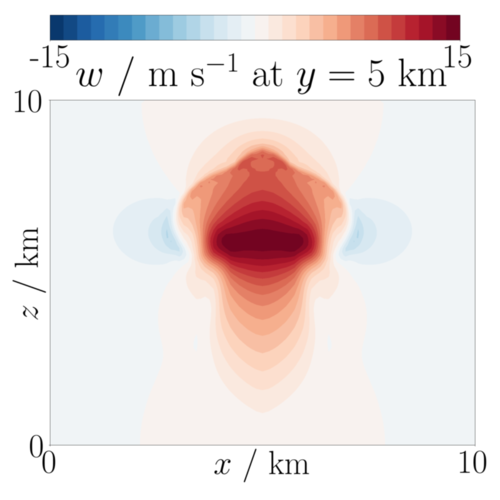}}
\caption[3D saturated rising bubble in $k=1$ configuration]{The outputted fields along $y=5$ km from the $k=1$ next-to-lowest order space simulation at $t=1000$ s of the three-dimensional simulation of a rising thermal in a saturated atmosphere. (Left) $\theta_e$ and (right) vertical velocity $w$ fields for a simulation using grid spacing $\Delta x=200$ m and a time step of $\Delta t=1$ s.
The $\theta_e$ field is contoured every $0.25$ K with the $320$ K contour omitted, whilst the contour spacing for the $w$ field is $1$ m s$^{-1}$.
Cross-sections are shown at $y=5$ km, with values plotted on the lower side of the plane.
As with the two-dimensional case, a second plume can be seen forming at the top of the primary plume.} \label{fig:3d bubble 2}
\end{figure}

We now demonstrate the use of our discretisation upon small-scale dynamics in three dimensions.
This test case is a three-dimensional version of the moist benchmark of \cite{bryan2002benchmark} that was described in Section \ref{sec:sat bubble}.
Rain and the effects of planetary rotation are not included. \\
\\
The domain is now periodic in the horizontal directions, with length, width and height $10$ km.
The background state set-up is the same as that in Section \ref{sec:sat bubble}, with $\theta_e=320$ K, $r_t=0.02$ kg kg$^{-1}$ and the pressure as $p=10^5$ Pa on the bottom surface.
Using the initialisation procedure that was outlined in Section \ref{sec:sat balance} generates the values of the prognostic variables such that the model is in hydrostatic balance and saturated with respect to water vapour. \\
\\
Using $x_c=y_c=5$ km and $z_c=2$ km, and defining
\begin{equation}
r = \sqrt{(x-x_c)^2+(y-y_c)^2+(z-z_c)^2},
\end{equation}
we apply the perturbation
\begin{equation}
\theta'_{vd}=\left\lbrace \begin{matrix}
\Delta\Theta\cos^2\left(\frac{\pi r}{2 r_c}\right), & r < r_c, \\
0, & \mathrm{otherwise},
\end{matrix}\right.
\end{equation}
with $\Delta \Theta = 1$ K and $r_c = 2$ km.
As in Section \ref{sec:sat bubble}, the perturbation is applied using the background field $\overline{\theta}_{vd}$
\begin{equation}
\theta_{vd} = \overline{\theta}_{vd}\left(1 + \frac{\theta'_{vd}}{300 \ \mathrm{K}}\right).
\end{equation}
The same routine as used in Section \ref{sec:sat bubble} is then applied to obtain the initial $\rho_d$, $r_v$ and $r_c$ fields, ensuring that the atmosphere is exactly saturated and that the initial pressure field is equal to the background pressure field. \\
\\
Cross-sections of the $\theta_e$ and vertical velocity $w$ fields at $t=1000$ s and $y=5$ km are shown in Figures \ref{fig:3d bubble 1} and \ref{fig:3d bubble 2}, for both the set-ups using $k=0$ lowest-order spaces (which had $\Delta x=\Delta y=\Delta z=100$ m) and the $k=1$ spaces (which had $\Delta x=\Delta y=\Delta z=200$ m).
Both simulations had $\Delta t=1$ s.
As in Section \ref{sec:sat bubble}, we see a secondary plume beginning to form at the top of the rising thermal in the $k=1$ case.

\subsection{Moist baroclinic wave}
\begin{figure}[h!]
\centering
{\includegraphics[width=0.4\textwidth]{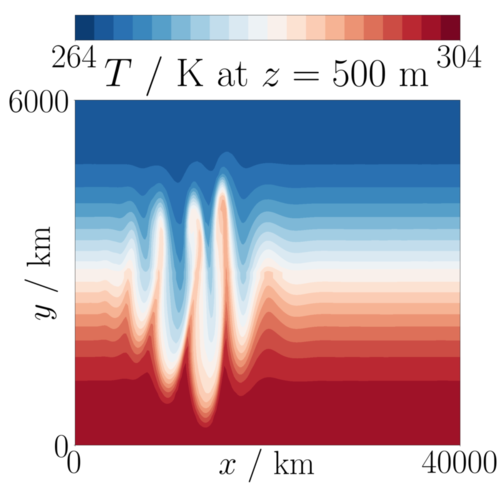}}~~
{\includegraphics[width=0.4\textwidth]{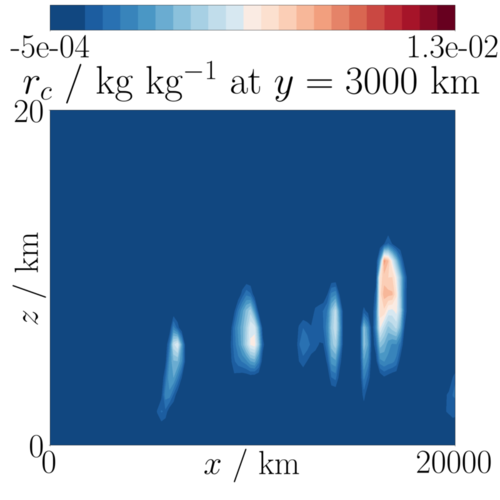}} \\
{\includegraphics[width=0.4\textwidth]{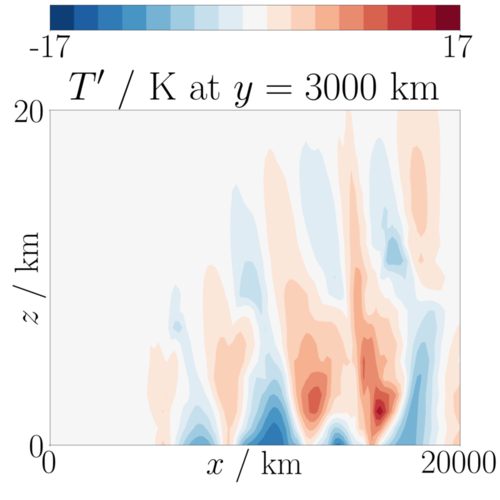}}~~
{\includegraphics[width=0.4\textwidth]{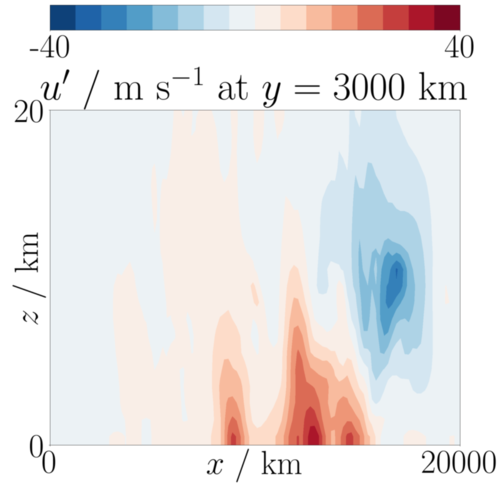}}
\caption[Moist baroclinic wave field outputs for $k=0$ configuration]{Cross-sections of fields at $t=12$ days from the moist baroclinic wave test case using the $k=0$ spaces.
The grid sizes used were $\Delta x=\Delta y=200$ km  and $\Delta z=1$ km, with $\Delta t=300$ s.
Shown are (top left) the $T$ field on $z=500$ m contoured every $2$ K, (top right) the $r_c$ field on $y=3000$ km contoured every $5\times 10^{-4}$ kg kg$^{-1}$ with the zero contour omitted, (bottom left)  the perturbed temperature field on $y=3000$ km with contours every $2$ K and (bottom right) the perturbed zonal wind $u'$, with contours spaced by $5$ m s$^{-1}$.
Values shown on the $y=3000$ km plane are computed from the lower side of the plane.} \label{fig:moist baroclinic 1}
\end{figure}

\begin{figure}[h!]
\centering
{\includegraphics[width=0.4\textwidth]{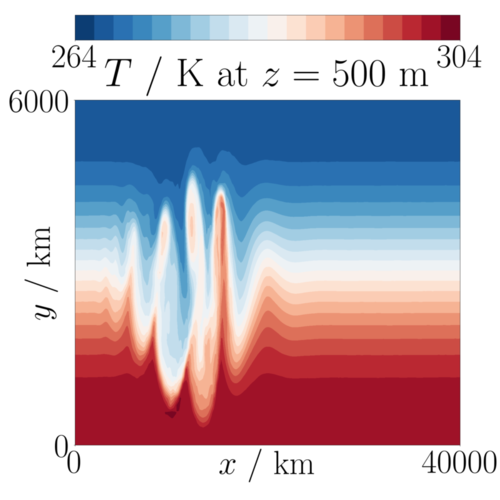}}~~
{\includegraphics[width=0.4\textwidth]{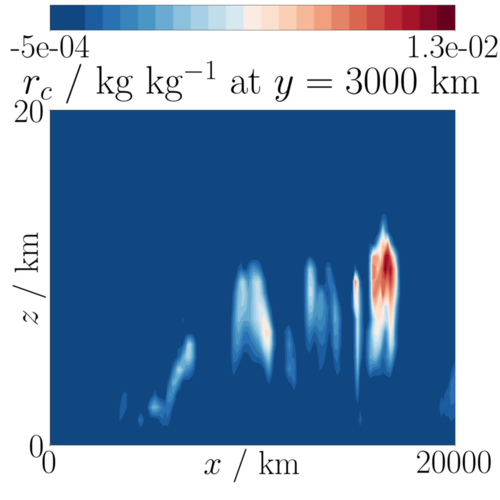}} \\
{\includegraphics[width=0.4\textwidth]{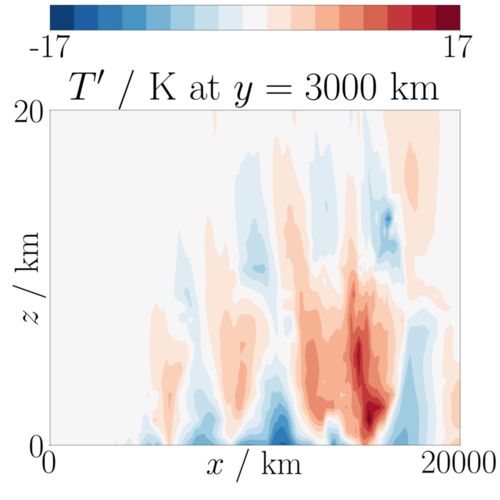}}~~
{\includegraphics[width=0.4\textwidth]{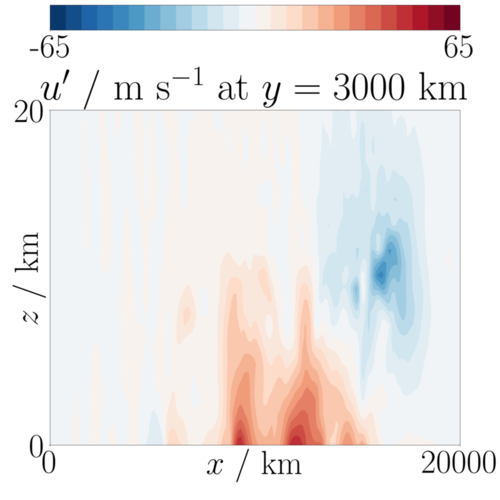}}
\caption[Moist baroclinic wave field outputs for $k=1$ configuration]{Cross-sections of fields at $t=12$ days from the moist baroclinic wave test case using the $k=1$ spaces.
The grid sizes used were $\Delta x=\Delta y=250$ km and $\Delta z=1.5$ km, with $\Delta t=300$ s.
Plots shown are the same as in Figure \ref{fig:moist baroclinic 1}, but note the different scale in the zonal wind perturbation plot (bottom right).} \label{fig:moist baroclinic 2}
\end{figure}

The final test case that we present is the moist baroclinic wave outlined in \cite{ullrich2015analytical}.
This is the only test featuring the Coriolis force, although rain is again neglected.
A three dimensional channel of length $L=40000$ km in the $x$-direction (in which the domain is periodic), width $W=6000$ km in the $y$-direction and height $H=30$ km in the $z$-direction.
The walls at the $y$ and $z$ boundaries of the domain are rigid, with no flow through them. \\
\\
This test case uses initial conditions that are analytically in thermal wind balance.
In \cite{ullrich2015analytical}, the vertical coordinate used is a pressure coordinate $\eta$.
This is used to define the background zonal wind $u$, the geopotential $\Phi$, the virtual temperature $T_v$ and the specific humidity $q$ according to the following equations:
\begin{subequations}
\begin{align}
u & = -u_0\sin^2\left(\frac{\pi y}{W}\right)\ln \eta \exp\left[-\left(\frac{\ln \eta}{b}\right)^2\right], \\
\Phi & = \frac{T_0 g}{\Gamma}\left(1- \eta^\frac{R_d \Gamma}{g}\right) + \frac{f_0 u_0}{2} \left[ y - \frac{W}{2}-\frac{W}{2\pi}\sin\left(\frac{2\pi y}{W}\right)\right]\ln \eta \exp \left[-\left(\frac{\ln \eta}{b}\right)^2\right], \\
T_v & = T_0\eta^\frac{R_d\Gamma}{g} + \frac{f_0 u_0}{2R_d} \left[ y - \frac{W}{2}-\frac{W}{2\pi}\sin\left(\frac{2\pi y}{W}\right)\right]\left[\frac{2}{b^2}(\ln \eta)^2-1 \right] \exp \left[-\left(\frac{\ln \eta}{b}\right)^2\right], \\
q & = \frac{q_0}{2}\exp\left[-\left(\frac{y}{\Delta y_w}\right)^4\right] \left \lbrace \begin{matrix}
1 + \cos\left[\frac{\pi(1-\eta)}{1-\eta_w}\right] & \eta \geq \eta_w, \\
0 & \mathrm{otherwise}.
\end{matrix} \right.
\end{align}
\end{subequations}
The constants take the values $\Gamma=0.005$ K m$^{-1}$, $f_0=2.00\times 10^{-6}$ s$^{-1}$, $a=6.37\times 10^6$ m, $T_0=288$ K, $u_0=35$ m s$^{-1}$, $b=2$, $\Delta y_w=3.2\times 10^6$ m and $\eta_w=0.3$. 
We use the slightly lower value of $q_0=0.016$ than \cite{ullrich2015analytical} to prevent our model being initially too close to saturation. \\
\\
This initial condition is converted to our prognostic variables using the Newton iteration procedure suggested in \cite{ullrich2015analytical} to find $\eta$.
Using the requirement that $\Phi=gz$ and taking $\eta\in V_\theta$, the procedure uses:
\begin{equation}
\eta^{n+1} = \eta^n - \frac{\Phi(\eta^n)-gz}{T_v(\eta^n)-R_d/\eta^n}.
\end{equation}
This $\eta$ field is then used to compute $u$, $T_v$ and $q$.
In order to convert $T_v$ and $q$ into $\theta_{vd}$ and $r_v$, we use the definitions of $\theta_{vd}$, $\eta=p/p_s$ for $p_s=10^5$ Pa and the equations
\begin{equation}
T = \frac{T_v}{1 + q(R_v/R_d - 1)}
\end{equation}
and
\begin{equation}
r_v = \frac{q}{1-q}.
\end{equation}
These $\theta_{vd}$ and $r_v$ fields are used to compute the $\rho_d$ field that provides hydrostatic balance for the background state using the procedure outlined in Section \ref{sec:hydrostatic balance}. \\
\\
The perturbation is added to the $x$-component of the background velocity, which we will denote as $\overline{u}$, so that $u=\overline{u}+u'$.
With $x_c=2\times10^6$ m, $y_c=2.5\times 10^6$, $L_p=6\times 10^5$, $u_p=1$ m s$^{-1}$ we use
\begin{equation}
r = \sqrt{(x-x_c)^2 + (y-y_c)^2},
\end{equation}
to get
\begin{equation}
u' = u_p \exp \left[-\left(\frac{r}{L_p}\right)^2 \right].
\end{equation}
Figures \ref{fig:moist baroclinic 1} and \ref{fig:moist baroclinic 2}  show fields from this test case at $t=12$ days for the $k=0$ and $k=1$ spaces respectively.
For the $k=0$ configuration, $\Delta x=\Delta y=200$ km  and $\Delta z=1$ km, whilst for $k=1$ these were $\Delta x=\Delta y=250$ km  and $\Delta z=1.5$ km.
For both simulations, $\Delta t=300$ s.
In the $k=1$ simulation, the baroclinic wave becomes much stronger than in the $k=0$ simulation.
As the wave develops, the maxima and minima in the temperatures in the $k=1$ case are higher than for the $k=0$ lowest-order spaces.
When these minima coincide with regions close to water vapour saturation, more condensation occurs.
This releases latent heat and strengthens the baroclinic wave, reinforcing the behaviour.


\section{Conclusions}
We have presented a discretisation of the moist compressible Euler equations that uses a compatible finite element framework.
Building upon the work of \cite{cotter2012mixed}, \cite{natale2016compatible}, \cite{yamazaki2017vertical} and \cite{bendall2019recovered}, we have presented two sets of compatible spaces that can be used.
The model configurations for each of the sets of spaces has been described, detailing the discrete equations that are solved and the transport schemes used, including slope limiters that can be used with moisture variables.
A discussion is given of the parametrisation of the moist processes and how these are coupled to the description of the resolved flow.
There is also a presentation of how a discrete hydrostatic balance can be set-up for the model under moist conditions.
The performance of the model is displayed through the presentation of five test cases, including two new ones which have introduced. \\
\\
Future work planned includes testing the choices made in these models and comparing the qualities of the model using each of the different sets of spaces.
More physics parametrisations will be added and the procedure used in the linear solve stage of the model will be thoroughly detailed along with its performance.

\section*{Acknowledgements}
TMB and THG were supported by the EPSRC Mathematics of Planet Earth Centre for Doctoral Training at Imperial College London and the University of Reading, with grant number EP/L016613/1.
\bibliographystyle{ieeetr}
\bibliography{moisture}
\section*{Appendix}
\textbf{Notation, thermodynamic variables and constants}
\begin{table}[h!]
\begin{tabular}{l|c|l}
Specific heat capacity of dry air at const. volume &
$c_{vd}$ &
717 J kg$^{-1}$ K$^{-1}$ \\
Specific heat capacity of dry air at const. pressure &
$c_{pd}$ &
1004.5 J kg$^{-1}$ K$^{-1}$ \\
Specific heat capacity of water vapour at const. volume &
$c_{vv}$ &
1424 J kg$^{-1}$ K$^{-1}$ \\
Specific heat capacity of water vapour at const. pressure &
$c_{pv}$ &
1885 J kg$^{-1}$ K$^{-1}$ \\
Specific heat capacity of liquid water at const. pressure &
$c_{pl}$ &
4186 J kg$^{-1}$ K$^{-1}$ \\
Specific heat capacity of moist air at const. volume &
$c_{vml}$ &
$c_{vd} + r_vc_{vv} + (r_c+r_r)c_{pl}$ \\
Specific heat capacity of moist air at const. pressure &
$c_{pml}$ &
$c_{pd} + r_vc_{pv} + (r_c+r_r)c_{pl}$ \\
Specific gas constant for dry air &
$R_d$ &
287 J kg$^{-1}$ K$^{-1}$ \\
Specific gas constant for water vapour &
$R_v$ &
461 J kg$^{-1}$ K$^{-1}$ \\
Specific gas constant for moist air &
$R_m$ &
$R_d+r_vR_v$ \\
Reference latent heat of vaporisation of water at $T_R$ &
$L_{vR}$ &
2.5$\times 10^6$ J kg$^{-1}$  \\
Latent heat of vaporisation of water &
$L_v$ &
$L_{vR} - (c_{pl}-c_{pv})(T-T_R)$ \\
Reference temperature &
$T_R$ &
273.15 K \\
Reference pressure &
$p_R$ &
$10^5$ Pa \\
Constant in Tetens' formula &
$C^\mathrm{sat}_0$ &
610.9 Pa \\
Constant in Tetens' formula &
$C^\mathrm{sat}_1$ &
-17.27 \\
Constant in Tetens' formula &
$C^\mathrm{sat}_2$ &
35.86 K \\
Ratio of $R_d$ to $R_v$ &
$\epsilon$ &
0.623 \\
Ratio of $R_d$ to $c_{pd}$ &
$\kappa$ &
2/7
\end{tabular}
\end{table}

\end{document}